\documentclass{siamart171218}
\usepackage[utf8]{inputenc}
\usepackage{amssymb}
\usepackage{enumitem}

\numberwithin{theorem}{section}
\newsiamremark{remark}{Remark}
\newsiamremark{example}{Example}

\def\theequation{\arabic{section}.\arabic{equation}}
\numberwithin{equation}{section}
\newcommand{\fd}[1]{\mathcal{D}^{#1}_t}
\newcommand{\dfd}[1]{\mathcal{D}^{#1}_\tau}
\newcommand{\diff}{\triangledown_\tau}
\newcommand{\Ass}[1]{\textbf{\upshape A#1}}
\newcommand{\Mss}[1]{\textbf{\upshape M#1}}
\newcommand{\R}{\mathbb{R}}
\newcommand{\iprod}[1]{\langle#1\rangle}
\newcommand{\taumax}{\tau}
\newcommand{\rhomax}{\rho}
\newcommand{\Err}[1]{\widetilde{\Pi_{#1}}}
\newcommand{\Gloc}{G_{\mathrm{loc}}}
\newcommand{\Ghis}{G_{\mathrm{his}}}

\newcommand{\vecx}{\boldsymbol{x}}
\newcommand{\Lapl}{\triangle}
\newcommand{\defeq}{:=}
\newcommand{\gammaopt}{\gamma_{\textrm{opt}}}

\newcommand{\zd}{\,\mathrm{d}}
\newcommand{\abs}[1]{\left|#1\right|}
\newcommand{\absb}[1]{\big|#1\big|}
\newcommand{\absB}[1]{\Big|#1\Big|}

\newcommand{\bra}[1]{\left(#1\right)}
\newcommand{\brab}[1]{\big(#1\big)}
\newcommand{\braB}[1]{\Big(#1\Big)}
\newcommand{\brat}[1]{(#1)}
\newcommand{\kbra}[1]{\left[#1\right]}

\newcommand{\myinnerb}[1]{\big\langle#1\big\rangle}

\newcommand{\mynorm}[1]{\left\|#1\right\|}
\newcommand{\mynormb}[1]{\big\|#1\big\|}

\newcommand{\bl}[1]{{\color{blue}#1}}

\newcommand{\TheTitle}{A second-order scheme with nonuniform time steps for a linear reaction-subdiffusion problem}

\title{{\TheTitle}\thanks{Submitted to the editors DATE.
\funding{This work was funded by NSFC grants 11771035, 91430216, U1530401; a grant 1008-56SYAH18037
from NUAA Scientific Research Starting Fund of Introduced Talent and a grant DRA2015518
from 333 High-level Personal Training Project of Jiangsu Provipnce;
Australian Research Council grant DP140101193.}}}
\author{
Hong-lin Liao\thanks{Department of Mathematics,
Nanjing University of Aeronautics and Astronautics,
Nanjing, 211106, P. R. China.
(\email{liaohl@csrc.ac.cn,liaohl@nuaa.edu.cn}).}
\and
William McLean\thanks{School of Mathematics and Statistics,
University of New South Wales, Sydney 2052, Australia.
(\email{w.mclean@unsw.edu.au}).}
\and
Jiwei Zhang\thanks{School of Mathematics and Statistics, and Hubei Key Laboratory of Computational Science, Wuhan University, Wuhan 430072, P. R. China. (\email{jiweizhang@whu.edu.cn}).}}

\begin{document}
\maketitle
\begin{abstract}
It is reasonable to assume that a discrete convolution structure dominates
the local truncation error of any numerical Caputo formula because the
fractional time derivative and its discrete approximation have the same
convolutional form. We suggest an error convolution structure (ECS) analysis
for a class of interpolation-type approximations to the Caputo fractional
derivative. Our assumptions permit the use of adaptive time steps, such as is
appropriate for accurately resolving the initial singularity of the solution
and also certain complex behavior away from the initial time. The ECS analysis
of numerical approximations has two advantages: (i) to localize (and
simplify) the analysis of the approximation error of a discrete convolution
formula on general nonuniform time grids; and (ii) to reveal the error
distribution information in the long-time integration via the global
consistency error. The core result in this paper is an ECS bound and a global
consistency analysis of the nonuniform Alikhanov approximation, which is
constructed at an offset point by using linear and quadratic polynomial
interpolation.  Using this result, we derive a sharp $L^2$-norm error estimate
of a second-order Crank-Nicolson-like scheme for linear reaction-subdiffusion
problems.  An example is presented to show the sharpness of our analysis.
\end{abstract}


\section{Introduction}\label{sec: introduction}
The time-fractional diffusion equation provides a valuable tool for modeling
complex systems such as glassy and disordered media~\cite{Hilfer:2000}. This
paper builds on our recent
results~\cite{LiaoLiZhang:2018,LiaoMcLeanZhang:2019,LiaoYanZhang:2019}
for the nonuniform mesh technique applied to the time discretization of the
following reaction-subdiffusion problem in a bounded domain $\Omega\subset\R^d$
($d=1$, $2$, $3$),
\begin{equation}\label{eq: IBVP}
\begin{aligned}
\fd{\alpha}u-\Lapl u&=\kappa u+f(\vecx,t)&
&\text{for $\vecx\in\Omega$ and $0<t<T$,}\\
u&=u_0(\vecx)&&\text{for $\vecx\in\Omega$ when $t=0$,}
\end{aligned}
\end{equation}
subject to the homogeneous Dirichlet boundary condition~$u=0$
on~$\partial\Omega$. Here, the reaction coefficient~$\kappa$ is a real constant,
and $\fd{\alpha}={}_{~0}^{C}\mathcal{D}_t^{\alpha}$ denotes the Caputo
fractional derivative of order~$\alpha$ ($0<\alpha<1$) with respect to time~$t$,
that is,
\[
(\fd{\alpha}v)(t)\defeq\int_0^t\omega_{1-\alpha}(t-s)v'(s)\zd{s}
	\quad\text{for $t>0$,}
	\quad\text{where~~ $\omega_\beta(t)\defeq t^{\beta-1}/{\Gamma(\beta)}$.}
\]

\subsection{Initial singularity and the nonuniform time meshes technique}
In developing numerical methods for solving the subdiffusion
problem~\eqref{eq: IBVP}, an important issue to be considered is that the
solution~$u$ is typically less regular than in the case of a classical
parabolic PDE (as the limiting case~$\alpha\to1$). Sakamoto and
Yamamoto~\cite{SakamotoYamamoto2011} showed that if the initial
data~$u^0\in H^2(\Omega)\cap H_0^1(\Omega)$, then the unique solution
$u\in C\bra{[0,T];H^2(\Omega)\cap H_0^1(\Omega)}$, with
$\fd{\alpha}u\in C\bra{[0,T];L^2(\Omega)}$~and $\partial_t u \in L^2(\Omega)$.
However, $\|\partial_t u(t)\|_{L^2(\Omega)} \le C_ut^{\alpha-1}$
for~$0<t\leq T$, where the constant~$C_u>0$ is independent of~$t$ but may depend
on~$T$. In fact, $u$ can only be a smooth function of~$t$ if the initial data
and source term satisfy some restrictive compatibility
conditions~\cite{Stynes2016}.

The focus of this paper is on a second-order time discretization
of~\eqref{eq: IBVP}.  The spatial discretization is of less interest: we apply
the standard Galerkin finite element method based on the weak form of the
fractional PDE,
\[
\iprod{\fd{\alpha}u,v}+\iprod{\nabla u,\nabla v}
	=\kappa\iprod{u,v}+\iprod{f(t),v}
	\quad\text{for all $v\in H^1_0(\Omega)$ and for $0<t\le T$,}
\]
where $\iprod{u,v}$ denotes the usual inner product in~$L_2(\Omega)$.
Thus, we construct the usual space of continuous, piecewise-linear functions
with respect to a partition of~$\Omega$ into subintervals (in 1D), triangles (in
2D) or tetrahedrons (in 3D) with the maximum diameter~$h$, and let $X_h$ denote
the subspace of functions satisfying the homogeneous Dirichlet boundary
condition. In the usual way, the (semidiscrete) Galerkin finite element
solution~$u_h:[0,T]\to X_h$ is then defined by requiring that
\[
\iprod{\fd{\alpha}u_h,\chi}+\iprod{\nabla u_h,\nabla\chi}
	=\kappa\iprod{u_h,\chi}+\iprod{f(t),\chi}
	\quad\text{for all $\chi\in X_h$ and for $0<t\le T$,}
\]
with $u_h(0)=u_{0h}\approx u_0$ for a suitable $u_{0h}\in X_h$.

Consider (generally nonuniform) time levels $0=t_0<t_1<t_2<\cdots<t_N=T$
and define a fractional time
level~$t_{n-\theta}\defeq\theta t_{n-1}+(1-\theta)t_n$ for an off-set parameter
$\theta\in[0,1/2)$. We denote the $k$th time-step size
by~$\tau_k\defeq t_k-t_{k-1}$ for~$1\le k\le N$, and the maximum step size by
$\taumax\defeq\max_{1\le k\le N}\tau_k$. We also define the local step-size
ratios
\[
\rho_k\defeq\frac{\tau_k}{\tau_{k+1}}\quad\text{for $1\le k\le N-1$,}\quad
\text{and put}\quad \rho\defeq\max_{1\le k\le N-1}\rho_k.
\]

For any time sequence $\brat{v^k}_{k=0}^N$, define the backward
difference~$\diff v^k:=v^k-v^{k-1}$ and the interpolated
value~$v^{n-\theta}\defeq\theta v^{n-1}+(1-\theta)v^n$. We consider a numerical
Caputo formula approximating $(\fd{\alpha}v)(t_{n-\theta})$ of the form
\begin{align}\label{eq: D tau}
(\dfd{\alpha}v)^{n-\theta}:=\sum_{k=1}^nA_{n-k}^{(n)}\diff v^k\approx
\sum_{k=1}^{n}\int_{t_{k-1}}^{\min\{t_k,t_{n-\theta}\}}\omega_{1-\alpha}(t_{n-\theta}-s)v'(s)\zd s\bl{,}
\end{align}
for appropriate discrete convolution kernels~$A_{n-k}^{(n)}$. Our fully-discrete
numerical solution, $u^n_h(\vecx)\approx u(\vecx,t_n)$
for~$\vecx\in\Omega$, is then defined by a time-stepping scheme: we require
that $u^n_h\in X_h$ satisfies
\begin{equation}\label{eq: discrete IBVP}
\iprod{(\dfd{\alpha}u_h)^{n-\theta},\chi}
	+\iprod{\nabla u_h^{n-\theta},\nabla\chi}
        =\kappa\iprod{u_h^{n-\theta},\chi}+\iprod{f(t_{n-\theta}),\chi}
\end{equation}
for all $\chi\in X_h$~and for~$1\le n\le N$, with~$u^0_h=u_{0h}$.

In the literature, several high-order numerical Caputo formulas have a discrete
convolution form like \eqref{eq: D tau}, such as the L1-2 schemes
\cite{GaoSunZhang:2014, LiaoLyuVongZhao:2016, LvXu:2016} and the
L2-1${}_\sigma$ formula \cite{Alikhanov2015,LiaoZhaoTeng:2016} that applied
the piecewise quadratic polynomial interpolation. They achieve second-order
temporal accuracy for sufficiently smooth solutions when applied to time
approximation of the pure subdiffusion equation~\eqref{eq: IBVP} with
$\kappa=0$. This article considers the L2-1${}_\sigma$ formula of
Alikhanov~\cite{Alikhanov2015}, which employs a quadratic interpolant
in each subinterval $[t_{k-1},t_{k}]$ for $1\leq k\leq n-1$, and a linear
interpolant in the final subinterval $[t_{n-1},t_{n-\theta}]$. The offset
parameter is chosen as~$\theta=\alpha/2$ (in our notation). As described below,
in the limit as~$\alpha\to1$, this scheme reduces to the well-known
Crank--Nicolson method ($\theta\rightarrow1/2$) for the classical diffusion equation.  We therefore refer to the time-stepping
scheme~\eqref{eq: discrete IBVP} as a \emph{fractional Crank--Nicolson} method.

In the special case of uniform time steps~$\tau_n=\tau$,
the discrete kernels $A^{(n)}_{n-k}=A_{n-k}$ depend only on the
difference~$n-k$, and were shown to be positive and monotonically decreasing,
leading to a proof that the resulting fractional Crank--Nicolson scheme is
stable and convergent of order~$O(\tau^2+h^2)$ in the $L_2$-norm assuming that
the solution~$u$ is sufficiently smooth~\cite{Alikhanov2015}. However, as
remarked above, in practice the time derivative~$\partial_tu$ typically behaves
like~$O(t^{\alpha-1})$ as~$t\to0$ \cite{SakamotoYamamoto2011,Stynes2016}, and so
this error bound breaks down.

In resolving a fixed singularity at~$t=0$ of the type described above, a simple
but useful technique to recover an optimal convergence order is to employ a
smoothly graded mesh $t_k=T(k/N)^{\gamma}$, where the grading
parameter~$\gamma\ge1$ is adapted to the strength of the singularity. The larger
the value of~$\gamma$ the more strongly the mesh points are concentrated
near~$t=0$. Actually, such meshes have long been used in the numerical solution
of Fredholm~\cite{Graham1982} and Volterra~\cite{Brunner1985} integral
equations, and their use for time-fractional PDEs is now well established
\cite{LiaoLiZhang:2018,LiaoLyuVongZhao:2016,McLeanMustapha2007,StynesEtAl2017}.
By using such a nonuniform mesh we will restore the second-order convergence in
time of the fractional Crank--Nicolson scheme in~\cite{Alikhanov2015} when the
solution is not smooth at~$t=0$. This idea was tested recently
in~\cite{LiaoZhaoTeng:2016} to resolve the initial singularity for the
subdiffusion problem, corresponding to~$\kappa=0$ in~\eqref{eq: IBVP}.
However, this is only a part of our story.

We will establish the stability and convergence theory for the fractional
Crank--Nicolson scheme on a wider class of nonuniform time meshes, not just the
standard graded mesh described above.  In this way, the theory could be
applied in advanced studies on adaptive time grids required to resolve certain
complex behavior (such as physical oscillations, blowup and so on) in nonlinear
time-fractional PDEs. These goals are natural, at least for linear
reaction-subdiffusion equations, since the backward Euler and Crank--Nicolson
schemes for the linear parabolic equation are stable and convergent (provided
$\tau\rightarrow0$) on arbitrary nonuniform grids with $\rho=O(1)$.

We refer the reader to other high-order time approximations in
\cite{JinLiZhou:2017,ZengEtAl2015} and the recent survey
paper~\cite{FordYan:2017}, which describes some useful approaches other than the
nonuniform grids technique to achieve second-order accuracy in time.

\subsection{Error convolution structure (ECS) analysis and a new problem}\label{subsec: ECS}
Generally, our goals are theoretically challenging because the numerical Caputo formula
always has a form of discrete convolutional summation \eqref{eq: D tau}.
Actually, the consistency analysis over the whole time interval $[t_0,t_{n-\theta}]$ becomes too cumbersome to implement in practice when there has not enough grid information.
To evade this difficulty, we propose an
\emph{error convolution structure} (ECS) analysis which begins by recasting
the discrete Caputo formula \eqref{eq: D tau} as
\begin{align}
(\dfd{\alpha}v)^{n-\theta}=A_{0}^{(n)}v^n
-\sum_{k=1}^{n-1}\brab{A_{n-k-1}^{(n)}-A_{n-k}^{(n)}}v^k-A_{n-1}^{(n)}v^{0}\,.
\end{align}
Consider the local truncation error
$\Upsilon^{n-\theta}:=(\fd{\alpha}v)(t_{n-\theta})-(\dfd{\alpha}v)^{n-\theta}$.
Given the construction of~$(\dfd{\alpha}v)^{n-\theta}$ via local interpolation
of~$v$, and provided the discrete convolution kernels $A_{k}^{(n)}$ are
decreasing, it is reasonable to conjecture that a discrete convolution structure
dominates the local truncation error:
\[
\text{(\textbf{ECS} hypothesis)}\qquad|\Upsilon^{n-\theta}|
    \leq A_{0}^{(n)}\Gloc^n
+\sum_{k=1}^{n-1}\brab{A_{n-k-1}^{(n)}-A_{n-k}^{(n)}}\Ghis^k.
\]
Here, $\Gloc^n$ arises from the interpolation error on the \emph{local}
subinterval $[t_{n-1},t_{n-\theta}]$ whereas the $\Ghis^k$ ($1\leq k\leq n-1$)
arise from the interpolation errors over the \emph{history} $[t_{0},t_{n-1}]$.
Obviously, this \textbf{ECS} hypothesis localizes the consistency analysis of discrete Caputo formulas,
and makes it possible to analyze the numerical approximations on a general class of nonuniform time grids.

Always, there is a loss of accuracy for $\Upsilon^{n-\theta}$ due to the initial
singularity of solution. Actually, $\Upsilon^{1-\theta}=O(1)$ holds on any
mesh and a superconvergence analysis should be required. For example, Stynes et
al.~\cite[Lemma 5.2]{Stynes2016} showed that, on a graded mesh, the truncation
error of the well-known L1 formula ($\theta=0$) behaves like
$\Upsilon^{n}=O\brab{n^{-\min\{2-\alpha,\gamma\alpha\}}}$. Building on the
ideas first introduced by Liao et al.~\cite[Section 3]{LiaoLiZhang:2018}, we
will prove a sharp error estimate via a fractional discrete Gronwall inequality
(Theorem \ref{thm: Gronwall}) that provides a \emph{global consistency error} in
the form
\begin{equation}\label{eq: globalConsistencyError-definition}
\mathcal{E}_{\mathrm{glob}}^n:=\sum_{k=1}^nP^{(n)}_{n-k}|\Upsilon^{k-\theta}|,
\quad\text{for $1\le k\le n\le N$,}
\end{equation}
where the \emph{complementary discrete convolution kernels}~$P^{(n)}_{n-k}$
are chosen to enforce the identity
\begin{equation}\label{eq: complementaryKernel-equality}
\sum_{j=k}^nP^{(n)}_{n-j}A^{(j)}_{j-k}\equiv1\quad\text{for $1\le k\le n\le N$.}
\end{equation}
In fact, rearranging this identity yields a recursive formula (in effect, a
definition)
\begin{align}\label{discreteConvolutionKernel-RL}
P_{0}^{(n)}\defeq\frac{1}{A_0^{(n)}},\quad
P_{n-j}^{(n)}\defeq
\frac{1}{A_0^{(j)}}
\sum_{k=j+1}^{n}\brab{A_{k-j-1}^{(k)}-A_{k-j}^{(k)}}P_{n-k}^{(n)},
    \quad 1\leq j\leq n-1.
\end{align}
In our recent paper~\cite{LiaoMcLeanZhang:2019}, we showed that this approach
is not limited to the L1 rule, but applies to a general class of discrete
convolution kernels~$A^{(n)}_{n-k}$~satisfying the following two assumptions:
\begin{description}
\item[\Ass{1}.] There is a constant~$\pi_A>0$ such that
$$A^{(n)}_{n-k}\ge\frac1{\pi_A\tau_k}\int_{t_{k-1}}^{t_k}
\omega_{1-\alpha}(t_n-s)\zd{s}\quad \text{for $1\le k\le n\le N$};$$
\item[\Ass{2}.] The discrete kernels are monotone,
$A^{(n)}_{k-2}\ge A^{(n)}_{k-1}>0$ for~$2\leq k\leq n\leq N$.
\end{description}
In this case, the complementary kernels~$P^{(n)}_{n-k}$ in
\eqref{discreteConvolutionKernel-RL} are well-defined and non-negative, and
satisfy~\cite[Lemma~2.1]{LiaoMcLeanZhang:2019}
\begin{align}\label{eq: P bound}
\sum_{j=1}^nP^{(n)}_{n-j}\omega_{1+(m-1)\alpha}(t_j)\leq \pi_A\omega_{1+m\alpha}(t_n)\quad\text{for $m=0,1$ and $1\le n\le N$.}
\end{align}
From the \textbf{ECS} hypothesis, one can exchange the order of summation
to find
\begin{align*}
\mathcal{E}_{\mathrm{glob}}^n
    &\leq\sum_{k=1}^nP^{(n)}_{n-k}A_{0}^{(k)}\Gloc^k
    +\sum_{k=1}^nP^{(n)}_{n-k}\sum_{j=1}^{k-1}
        \brab{A_{k-j-1}^{(k)}-A_{k-j}^{(k)}} \Ghis^j\\
    &\leq\sum_{k=1}^nP^{(n)}_{n-k}A_{0}^{(k)}\Gloc^k
    +\sum_{j=1}^{n-1}\Ghis^j\sum_{k=j+1}^{n}
    \brab{A_{k-j-1}^{(k)}-A_{k-j}^{(k)}}P_{n-k}^{(n)},
\end{align*}
and then, by using the definition \eqref{discreteConvolutionKernel-RL}
directly, arrive at
\begin{align}\label{ieq: globalConsistencyBound}
\mathcal{E}_{\mathrm{glob}}^n\leq&\,\sum_{k=1}^nP^{(n)}_{n-k}A_{0}^{(k)}\Gloc^k
+\sum_{k=1}^{n-1}P^{(n)}_{n-k}A_{0}^{(k)}\Ghis^k\,.
\end{align}
Thus, by using the properties in \eqref{eq: complementaryKernel-equality}~and
\eqref{eq: P bound}, it is possible to obtain some useful error estimates on
a variety of nonuniform grids, not limited to~$t_k=T(k/N)^\gamma$.

Obviously, the first term on the right-hand side
of~\eqref{ieq: globalConsistencyBound} represents the total error contributions
from discretization errors over the $n$~current cells $[t_{k},t_{k-\theta}]$
($1\leq k\leq n$), whereas the second term represents those from discretization
errors over the $\frac{1}{2}n(n-1)$ small cells in the historic intervals
$[t_{0},t_{k-1}]$ ($2\leq k\leq n$). This observation is very interesting: the
local error in the current cell $[t_{n},t_{n-\theta}]$ and the historic errors
in the (long-time) interval $[t_0,t_{n-1}]$ make almost the same contribution,
in the sense of convolutional summation, to the global consistency error of the
discrete Caputo derivative. If some appropriate time grid is chosen to make $\Ghis^k=O\bra{\Gloc^k}$
according to the error equidistribution  principle,
the error bound \eqref{ieq: globalConsistencyBound} becomes
\begin{align*}
\mathcal{E}_{\mathrm{glob}}^n\lesssim\sum_{k=1}^nP^{(n)}_{n-k}A_{0}^{(k)}\Gloc^k\,.
\end{align*}
It suggests that the global approximation error of the numerical Caputo
formula \eqref{eq: D tau} depends mainly on the local error $\Gloc^k$. In this sense,
\emph{the error of numerical Caputo
formula  is ``local" despite its overtly nonlocal nature.}

\begin{table}[!ht]
\caption{Mesh restriction to stability for linear reaction-(sub)diffusion equations}
\label{tab: mesh condition}
\begin{center}
\renewcommand{\arraystretch}{1.1}
\begin{tabular}{c|c|c}
  & backward Euler-like & Crank-Nicolson-like \\
  \hline
diffusion ($\alpha\rightarrow1$)   & $\rho=O(1)$ & $\rho=O(1)$ \\
\hline
subdiffusion ($0<\alpha<1$) & $\rho=O(1)$ & ? \\
\hline
\end{tabular}
\end{center}
\end{table}

The \textbf{ECS} hypothesis plays a key role in our analysis. Actually, it has been used implicitly for the nonuniform L1 (fractional backward
Euler-type) method employing a linear interpolant
in each subinterval $[t_{k-1},t_{k}]$ for $1\leq k\leq n$.
That analysis \cite[(3.9) in Lemma 3.3]{LiaoLiZhang:2018}
showed that the \textbf{ECS} hypothesis
is valid for~$\rho=1$, or in other words provided $\tau_k\le\tau_{k+1}$ for
all~$k$. In a further study \cite{LiaoYanZhang:2019} on the two-level fast L1
scheme (which includes the original L1 scheme as a special case by setting the
SOE approximation error $\epsilon\equiv0$), the \textbf{ECS} hypothesis is
shown to be valid for any nonuniform mesh with~$\rho=O(1)$
\cite[Lemma 3.1]{LiaoYanZhang:2019}, that is,
\begin{align*}
|\Upsilon^{n}|\leq
a_{0}^{(n)}G^n+\sum_{k=1}^{n-1}\brab{a_{n-k-1}^{(n)}-a_{n-k}^{(n)}}G^k
\end{align*}
where the L1 kernels $a_{n-k}^{(n)}\defeq\frac{1}{\tau_k}\int_{t_{k-1}}^{t_k}\omega_{1-\alpha}(t_n-s)\zd{s}$ for $1\le k\leq n$,
and
$$G^k :=2\int_{t_{k-1}}^{t_k}\bra{t-t_{k-1}}\abs{v''(t)}\zd t\quad \text{for $1\le k\leq n$.}$$
This local step ratio restriction is the same as
that for the backward Euler scheme for a classical diffusion equation.
Considering Table~\ref{tab: mesh condition}, it is then natural to ask
an elementary problem: what restriction on~$\rho$ will suffice to ensure that the  fractional
Crank--Nicolson time-stepping scheme~\eqref{eq: discrete IBVP} is stable and
convergent? We address this problem in the condition \textbf{M1} below.

\subsection{The discrete fractional Gr\"{o}nwall inequality and our answer}

Our answer relies also on a discrete fractional Gr\"{o}nwall inequality suited
to general nonuniform time meshes, proved in our recent
paper~\cite[Theorem~3.1]{LiaoMcLeanZhang:2019} and
stated below (in a simplified form).  This result involves the aforementioned
complementary discrete convolution kernels~$P^{(n)}_{n-k}$, which are
well-defined thanks to our assumptions \Ass{1}--\Ass{2} on the discrete
convolution kernels~$A^{(n)}_{n-k}$ in the numerical Caputo
formula~\eqref{eq: D tau}.  The Mittag--Leffler function
$E_\alpha(z)\defeq\sum_{k=0}^\infty z^k/\Gamma(1+k\alpha)$ also appears.
\begin{theorem}\label{thm: Gronwall}
Let the criteria \Ass{1}--\Ass{2} hold, and the offset parameter $\theta\in[0,1)$.
Suppose that $\lambda>0$ is a constant independent of the time steps
and that the maximum time-step size
\[\taumax\le1/\sqrt[\alpha]{2\Gamma(2-\alpha)\pi_A\lambda}.\]
If the non-negative time sequences $(\xi^k)_{k=1}^N$ and $(v^k)_{k=0}^N$ satisfy
\begin{equation}\label{eq: first Gronwall}
\sum_{k=1}^nA^{(n)}_{n-k}\diff\bra{v^k}^2\le
	\lambda\bra{v^{n-\theta}}^2+v^{n-\theta}\xi^n\quad\text{for $1\le n\le N$,}
\end{equation}
then the solution $\brat{v^k}_{k=0}^N$ satisfies, for $1\le n\le N$,
\begin{align}
v^n&\le2E_\alpha\bigl(2\max(1,\rhomax)\pi_A \lambda t_n^\alpha\bigr)
	\braB{v^0+\max_{1\le k\le n}\sum_{j=1}^k P^{(k)}_{k-j}\xi^j
	}\label{eq: Gronwall conclusion}\\
&\le2E_\alpha\bigl(2\max(1,\rhomax)\pi_A \lambda t_n^\alpha\bigr)
	\braB{v^0+\pi_A\Gamma(1-\alpha)\max_{1\le j\le n}\{t_j^{\alpha}\xi^j\}}\,.\label{eq: Gronwall conclusion 1}
\end{align}
\end{theorem}

Thus, we need to complete the following three tasks:
\begin{description}
\item[Task 1.] Verify the assumptions \Ass{1}--\Ass{2} for the nonuniform
Alikhanov kernels~$A^{(n)}_{n-k}$ (defined in \cref{sec: Alikhanov})
so that we can use the complementary kernels~$P^{(n)}_{n-k}$ and apply the
fractional Gr\"{o}nwall inequality to establish the stability of the fully
discrete scheme \eqref{eq: discrete IBVP}.
\item[Task 2.] Verify the \textbf{ECS} hypothesis on nonuniform time meshes
and determine the corresponding expressions for $\Gloc^k$~and $\Ghis^k$ to
insert in the bound \eqref{ieq: globalConsistencyBound} for the global
consistency error $\mathcal{E}_{\mathrm{glob}}^n$.
\item[Task 3.] Establish a sharp error estimate in~$L^2$ for the fully discrete
scheme \eqref{eq: discrete IBVP} for the subdiffusion problem \eqref{eq: IBVP}
taking the initial singularity into account.
\end{description}

In more detail, we complete Task~1 in~\cref{sec: Alikhanov}. We describe the
fractional Crank--Nicolson scheme and the corresponding discrete Alikhanov
kernels~$A^{(n)}_{n-k}$, and show in \cref{thm: A1 A2} (The lengthy and technical proofs for these properties of the discrete
kernels~$A^{(n)}_{n-k}$ are postponed until~\cref{sec: monotone}) that the criteria
\Ass{1}--\Ass{2} hold given the following assumption on the mesh.

\begin{description}
\item[\Mss{1}.] The parameter $\theta=\alpha/2$, and the maximum time-step
ratio is $\rho=7/4$.
\end{description}

The special choice of~$\theta$ in \Mss{1} is needed in any case to achieve
second-order accuracy (see \cref{rem: theta}). At the end
of \cref{sec: Alikhanov}, the discrete fractional Gr\"{o}nwall inequality is applied to
establish stability for the time-stepping scheme \eqref{eq: discrete IBVP}.
Actually, by showing that $v^n=\|u^n_h\|$ satisfies \eqref{eq: first Gronwall},
the \emph{a priori} estimate with respect to initial and external
perturbations in the
forms~\eqref{eq: Gronwall conclusion}--\eqref{eq: Gronwall conclusion 1},
follows.

To verify the \textbf{ECS} hypothesis in Task~2 we make use of a proper lower
bound for $A^{(n)}_{n-k-1}-A^{(n)}_{n-k}$, already proved in part~(II)
of~\cref{thm: A1 A2} to ensure the criterion \Ass{2} directly. In the first
part of \cref{sec: truncation}, an interpolation error formula for quadratic
polynomials is derived in \cref{lem:quadraticInterpolationError}. Then we
complete Task~2 in \cref{thm: Upsilon global} by showing that the \textbf{ECS}
hypothesis and the bound \eqref{ieq: globalConsistencyBound} for the global
consistency error~$\mathcal{E}_{\mathrm{glob}}^n$ are valid under the
condition~\Mss{1}.

Task 3 is completed in the second part of \cref{sec: truncation}.
To make our analysis extendable (such as, for distributed-order
subdiffusion problems), we assume that there is
a constant $C_u>0$ such that the continuous solution $u$ satisfies
\begin{equation}\label{eq: sigma}
\begin{aligned}
\|u^{(l)}(t)\|_{H^2(\Omega)}\le C_u\brab{1+t^{\sigma-l}} \;\; \quad \text{for $l = 0,1,2,3$, and $0<t\le T$,}
\end{aligned}
\end{equation}
where $\sigma\in(0,1)\cup(1,2)$ is a regularity parameter. For
example~\cite{McLean2010,SakamotoYamamoto2011,StynesEtAl2017},
the assumption \eqref{eq: sigma} holds with $\sigma=\alpha$ for the subdiffusion problem
\eqref{eq: IBVP} if $f(\vecx,t)=0$ and $u_0\in H^1_0(\Omega)\cap H^2(\Omega)$.
To resolve such a solution $u$ efficiently, it is appropriate to choose the time
mesh in such a way that the following
condition~\cite{Brunner1985,McLeanMustapha:2015} holds.
\begin{itemize}
\item[\Mss{2}.] There is a constant $C_{\gamma}>0$ such that
$\tau_k\le C_{\gamma}\tau\min\{1,t_k^{1-1/\gamma}\}$
for~$1\le k\le N$, with $t_{k}\leq C_{\gamma}t_{k-1}$~and
$\tau_k/t_k \leq  C_{\gamma} \tau_{k-1}/t_{k-1}$ for~$2\le k\le N$.
\end{itemize}
Here, the parameter $\gamma\ge1$ controls the extent to which the time levels
are concentrated near $t=0$.  If the mesh is quasi-uniform, then
\Mss{2} holds with~$\gamma=1$.  As~$\gamma$
increases, the initial step sizes become smaller compared to the later ones.
A simple example of a family of meshes satisfying \Mss{2} is the graded
mesh~$t_k=T(k/N)^{\gamma}$.

When the offset parameter $\theta=0$ and \eqref{eq: D tau} is the nonuniform L1
method, our previous work~\cite[Theorem~3.1]{LiaoLiZhang:2018} proved the
following error bound for the fully discrete scheme~\eqref{eq: discrete IBVP},
\[
\mynormb{u(t_n)-u^n_h}\leq
 \frac{C_u}{\sigma(1-\alpha)}\tau^{\min\{\gamma\sigma,2-\alpha\}}+C_uh^2,\quad
1\leq n\leq N.
\]
In particular, the error is of order~$O(\tau^{2-\alpha}+h^2)$ if
$\gamma\geq(2-\alpha)/\sigma$. When $\theta=\alpha/2$ and
\eqref{eq: D tau} is the Alikhanov formula, \cref{thm: convergence}
establishes an error bound
\begin{equation}\label{eq: error bound}
\mynormb{u(t_n)-u^n_h}\le
    \frac{C_u}{\sigma(1-\alpha)}\tau^{\min\{\gamma\sigma,2\}}+C_uh^2,
    \quad 1\leq n\leq N,
\end{equation}
which is of order $O(\tau^{2}+h^2)$ if $\gamma\ge2/\sigma$.  Thus, in comparison
to the L1 scheme, the Alikhanov formula leads to a higher convergence rate with
respect to~$\tau$; however, both methods achieve only order
$O(\tau^{\sigma}+h^2)$ convergence on a uniform mesh. Numerical experiments in
\cref{sec: numerical} confirm that our error bound~\eqref{eq: error bound} is
sharp.

\section{Numerical Caputo formula and stability}\label{sec: Alikhanov}

Let $\Pi_{1,k}v$ denote the linear interpolant of a function~$v$ with
respect to the nodes $t_{k-1}$~and $t_k$, and let $\Pi_{2,k}v$ denote the
quadratic interpolant with respect to $t_{k-1}$, $t_k$~and $t_{k+1}$.
The corresponding interpolation errors are denoted by
\[
(\Err{p,k}v)(t)\defeq v(t)-\bra{\Pi_{p,k}v}(t)\quad
\text{for} \quad p\in\{1,2\}.
\]
Recalling that $\rho_k=\tau_k/\tau_{k+1}$, it is easy to find (for instance, by using
the Newton forms of the interpolating polynomials) that
\[
\bra{\Pi_{1,k}v}'(t)=\frac{\diff v^k}{\tau_k}
\;\;\text{and}\;\;
\bra{\Pi_{2,k}v}'(t)=\frac{\diff v^k}{\tau_k}
	+\frac{2(t-t_{k-1/2})}{\tau_k(\tau_k+\tau_{k+1})}
       \bra{\rho_k\diff v^{k+1}-\diff v^k}.
\]
Throughout this paper, we will always use the notation
\[
\varpi_n(t):=-\omega_{2-\alpha}(t_{n-\theta}-t)\le0
    \quad\text{for $0\leq t\leq t_{n-\theta}$}.
\]
If $0\leq t<t_{n-\theta}$, then
$\varpi_n'(t)=\omega_{1-\alpha}(t_{n-\theta}-t)>0$,
$\varpi_n''(t)=-\omega_{-\alpha}(t_{n-\theta}-t)>0$
and $\varpi_n'''(t)=\omega_{-\alpha-1}(t_{n-\theta}-t)>0$.
\subsection{Discrete Caputo formula}
The nonuniform Alikhanov approximation to the Caputo
derivative~$(\fd{\alpha}v)(t_{n-\theta})$ is defined by
\begin{align}\label{eq: L2-1_sigma}
(\dfd{\alpha}v)^{n-\theta}
	&\defeq\int_{t_{n-1}}^{t_{n-\theta}}\varpi_n'(s)
		\bra{\Pi_{1,n}v}'(s)\zd{s}	+\sum_{k=1}^{n-1}\int_{t_{k-1}}^{t_k}	\varpi_n'(s)\bra{\Pi_{2,k}v}'(s)\zd{s}\\
   &=a^{(n)}_0\diff v^n+\sum_{k=1}^{n-1}\braB{
	a^{(n)}_{n-k}\diff v^k+\rho_k b^{(n)}_{n-k}\diff v^{k+1}
		-b^{(n)}_{n-k}\diff v^k},\nonumber
\end{align}
where the discrete coefficients
$a_{n-k}^{(n)}$ and $b_{n-k}^{(n)}$ are defined by
\begin{align}
&a^{(n)}_{n-k}\defeq\frac{1}{\tau_k}\int_{t_{k-1}}^{\min\{t_k,t_{n-\theta}\}}\varpi_n'(s)\zd{s},\quad 1\le k\leq n;\label{eq: an}\\
&b^{(n)}_{n-k}\defeq\frac{2}{\tau_k(\tau_k+\tau_{k+1})}\int_{t_{k-1}}^{t_k}
	(s-t_{k-\frac12})\varpi_n'(s)\zd{s}, \quad 1\le k\leq n-1. \label{eq: bn}
\end{align}

When~$\theta=0$, the coefficients $a^{(n)}_{n-k}$ in \eqref{eq: an} are just the discrete convolution kernels
in the L1 formula~\cite{LiaoLiZhang:2018}. Notice that if $\alpha\rightarrow1$, then $\omega_{2-\alpha}(t)\to1$
whereas $\omega_{1-\alpha}(t)\to0$, uniformly for~$t$ in any compact
subinterval of the open half-line~$(0,\infty)$.  Thus,
$$a^{(n)}_0=\omega_{2-\alpha}((1-\theta)\tau_n)/\tau_n\to1/\tau_n$$
whereas $a^{(n)}_{n-k}\to0$~and $b^{(n)}_{n-k}\to0$ for $1\le k\le n-1$.
It follows that $(\dfd{\alpha}v)^{n-\theta}\to\diff{v^n}/\tau_k$
and $\theta=\alpha/{2}\to1/2$
so the scheme~\eqref{eq: discrete IBVP} tends to the
Crank--Nicolson method for a linear
parabolic equation. This is why we also call \eqref{eq: discrete IBVP} a fractional Crank--Nicolson time-stepping method.

Rearranging the terms in~\eqref{eq: L2-1_sigma}, we obtain
the compact form \eqref{eq: D tau}
where the discrete convolution kernels $A_{n-k}^{(n)}$ are defined as follows:
$A_0^{(1)}\defeq a_0^{(1)}$ if $n=1$ and, for $n\geq2$,
\begin{equation}\label{eq: weights}
A^{(n)}_{n-k}\defeq\begin{cases}
	a^{(n)}_0+\rho_{n-1}b^{(n)}_1,
	&\text{for $k=n$,}\\
	a^{(n)}_{n-k}+\rho_{k-1}b^{(n)}_{n-k+1}-b^{(n)}_{n-k},
	&\text{for $2\le k\le n-1$,}\\
	a^{(n)}_{n-1}-b^{(n)}_{n-1},
	&\text{for $k=1$.}
\end{cases}
\end{equation}

Before studying the kernels~$A^{(n)}_{n-k}$, we present two alternative formulas
for~$b_{n-k}^{(n)}$. Recall the integral form of error term for the
trapezoidal rule, which can be derived by the Taylor expansion with the
integral remainder.
Integration by parts yields the following identities.
\begin{lemma}\label{lem: integralFormula}
For any function $q\in C^2([t_{k-1},t_k])$,
\begin{align*}
\int_{t_{k-1}}^{t_k}(s-t_{k-1/2})q'(s)\zd{s}&=-\int_{t_{k-1}}^{t_k}\brab{\Err{1,k}q}(s)\zd{s}\\
	&=\frac{1}{2}\int_{t_{k-1}}^{t_k}(s-t_{k-1})(t_k-s)q''(s)\zd{s}.
\end{align*}
\end{lemma}
Taking $q:=\varpi_n$ in Lemma \ref{lem: integralFormula}, the definition \eqref{eq: bn} of $b_{n-k}^{(n)}$ gives
\begin{align}
b_{n-k}^{(n)}=&\,-2
\int_{t_{k-1}}^{t_{k}}\frac{\brab{\widetilde{\Pi_{1,k}}\varpi_n}(s)\zd{s}}{\tau_k(\tau_{k+1}+\tau_k)}
\label{eq: bn-anlter}\\
=&\, \int_{t_{k-1}}^{t_{k}}\frac{(t_k-s)(s-t_{k-1})}{\tau_k(\tau_{k+1}+\tau_k)}
\varpi_n''(s)\zd{s},\quad 1\leq k\leq n-1.\label{eq: bn-anlternative}
\end{align}

The following theorem gathers some useful properties of the
discrete kernels~$A^{(n)}_{n-k}$, but the rigorous proof is left to~\cref{sec: monotone}.
It should be noted here that this proof is quite different from the previous analysis
\cite{Alikhanov2015,GaoSunZhang:2014,LvXu:2016} for the discrete convolution kernels
in high-order numerical Caputo formulas with uniform time-steps.

\begin{theorem}\label{thm: A1 A2}
Let \Mss{1} hold and consider the discrete kernels defined
in~\eqref{eq: weights}.
\begin{itemize}
\item[(I)] The discrete kernels $A^{(n)}_{n-k}$ are bounded,
$$A_{0}^{(n)}\leq
\frac{24}{11\tau_n}\int_{t_{n-1}}^{t_{n}}\omega_{1-\alpha}(t_n-s)\zd s$$
and
 $$A_{n-k}^{(n)}\geq\frac4{11\tau_k}\int_{t_{k-1}}^{t_{k}}\omega_{1-\alpha}(t_n-s)\zd s,\quad 1\le k\le n;$$
\item[(II)] The discrete kernels $A^{(n)}_{n-k}$ are monotone,
\[A^{(n)}_{n-k-1}-A^{(n)}_{n-k}\ge
(1+\rho_k)b^{(n)}_{n-k}+\frac1{5\tau_k}\int_{t_{k-1}}^{t_k}\bra{t_k-s}\varpi_n''(s)\zd{s},
\quad 1\le k\le n-1;\]
\item[(III)] And the first kernel $A^{(n)}_{0}$ is appropriately larger than the second one,
 $$ \frac{1-2\theta}{1-\theta}A^{(n)}_{0}-A^{(n)}_{1}>0\quad\text{ for~$n\ge2$.}$$
\end{itemize}
\end{theorem}

The first part (I) implies that the criterion \Ass{1} holds with
$\pi_A=\frac{11}{4}$, the second part (II) ensures that the criterion \Ass{2} is
valid and the third part (III) is used to prove the following corollary. These
results allow us to apply \cref{thm: Gronwall} and establish the stability of
the time-stepping scheme~\eqref{eq: discrete IBVP}. Also, the second part (II)
establishes a stronger estimate used in obtaining an \textbf{ECS} bound for the
error analysis (see \cref{thm: Upsilon global}).
\begin{corollary}\label{cor: Dv v} Under the condition \Mss{1}, the discrete Caputo formula \eqref{eq: D tau}
with the discrete kernels \eqref{eq: weights} satisfies
$$\myinnerb{\bra{\dfd{\alpha}v}^{n-\theta},v^{n-\theta}}
    \geq\frac12\sum_{k=1}^nA^ {(n)}_{n-k}\diff\brab{\mynorm{v^k}^2}\quad
\text{for~$1\le n\le N$.}$$
\end{corollary}
\begin{proof}
The inequality is known to hold \cite[Lemma~4.1]{LiaoMcLeanZhang:2019}
provided \Ass{2} is satisfied and
$\theta^{(n)}\ge\theta$ for~$1\le n\le N$, where
\[\theta^{(1)}=\frac12\quad\text{and}\quad
\theta^{(n)}=\frac{A^{(n)}_0-A^{(n)}_1}{2A^{(n)}_0-A^{(n)}_1}
\quad\text{for $n\ge2$.}
\]
Obviously, \cref{thm: A1 A2} (II) ensures that \Ass{2} holds,
and the condition \Mss{1} leads to $\theta^{(1)}\ge\theta$.
From \cref{thm: A1 A2} (III), $\theta^{(n)}\geq\theta$ holds also for~$n\ge2$.
\end{proof}

\subsection{Unconditional stability}
By taking the $\chi=u_h^{n-\theta}$ in~\eqref{eq: discrete IBVP}, one has
\begin{equation}\label{eq: energy}
\myinnerb{\bra{\dfd{\alpha}u_h}^{n-\theta},u_h^{n-\theta}}\le
\kappa_{+}\mynormb{ u^{n-\theta}_h}^2+\myinnerb{f(t_{n-\theta}),u^{n-\theta}_h}
    \quad\text{for $1\le n\le N$},
\end{equation}
where $\kappa_{+}:=\max\{\kappa,0\}$ and the property
 {$\iprod{\nabla u_h^{n-\theta},\nabla u_h^{n-\theta}}\geq0$}
was used.  Therefore, applying the above \cref{cor: Dv v} along with the Cauchy--Schwarz
and triangle inequalities, one gets
\begin{multline*}
\sum_{k=1}^nA^{(n)}_{n-k}\diff\brab{\mynorm{u^k_h}^2}
    \le2\kappa_{+}\braB{(1-\theta)\mynormb{u^n_h}+\theta\mynormb{u^{n-1}_h}}^2\\
+2\braB{(1-\theta)\mynormb{u^n_h}+\theta\mynormb{u^{n-1}_h}}
    \mynormb{ {f(t_{n-\theta})}},\quad1\le n\le N,
\end{multline*}
which has the form of \eqref{eq: first Gronwall} with
\[
\lambda:=2\kappa_{+},\quad v^k:=\mynormb{u^k}\quad \text{and}\quad\xi^k:=2\mynormb{f^{k-\theta}}\quad\text{for $1\le k\le N.$}
\]
Note that \cref{thm: A1 A2} shows the criteria \Ass{1}--\Ass{2} of \cref{thm: Gronwall} are
satisfied with $\pi_A=11/4$, and the condition \Mss{1} gives~$\rho=7/4$. Therefore, applying \cref{thm: Gronwall}, we see that the
time-stepping method~\eqref{eq: discrete IBVP} is stable in the following sense.

\begin{theorem}\label{thm: stability}
If \Mss{1} holds with the maximum time step
$\taumax\le1/\sqrt[\alpha]{11\Gamma(2-\alpha)\kappa_{+}}$
(there is no limit to the maximum time step if $\kappa\le0$),
then the solution~$u^n_h$ of the time-stepping scheme~\eqref{eq: discrete IBVP}
is stable, that is,
\begin{align*}
\mynorm{u^n_h}&\le2E_\alpha\bigl(20\kappa_{+}t_n^\alpha\bigr)
\braB{\mynorm{ {u_{0h}}}+
	2\max_{1\le k\le n}\sum_{j=1}^kP^{(k)}_{k-j}
	\mynorm{ {f(t_{j-\theta})}}}\\
&\le2E_\alpha\bigl(20\kappa_{+}t_n^\alpha\bigr)
\braB{\mynorm{ {u_{0h}}}+
	6\Gamma(1-\alpha)\max_{1\le j\le n}\{t_j^{\alpha}
	\mynorm{ {f(t_{j-\theta})}}\}}\quad\text{for $1\le n\le N.$}
\end{align*}
\end{theorem}

\section{Global consistency error and convergence}\label{sec: truncation}

We now derive a representation for the consistency error of the discrete Caputo
derivative~\eqref{eq: D tau} with the discrete kernels in~\eqref{eq: weights}.
Fix a function~$v(t)$ and decompose the local consistency error into $n$~terms
corresponding to the $n$~subintervals, writing
\begin{equation}\label{eq: Upsilon}
\Upsilon^{n-\theta}:=(\fd{\alpha}v)(t_{n-\theta})-(\dfd{\alpha}v)^{n-\theta}
    =\sum_{k=1}^n\Upsilon^{n-\theta}_k,	\qquad1\le n\le N,
\end{equation}
where, recalling the notations~$\varpi_n(s)$, $(\Err{1,k}v)$ and $(\Err{2,k}v)$
from \cref{sec: Alikhanov},
\begin{align}
&\Upsilon^{n-\theta}_k\defeq\int_{t_{k-1}}^{t_k}
	\varpi_n'(s)\brab{\Err{2,k}v}'(s)\zd{s},\quad1\le k\le n-1\leq N-1,\label{eq: Upsilon n k}\\
&\Upsilon^{n-\theta}_n\defeq\int_{t_{n-1}}^{t_{n-\theta}}
	\varpi_n'(s)\brab{\Err{1,n}v}'(s)\zd{s},\quad 1\leq n\leq N.\label{eq: Upsilon n n}
\end{align}

Compared with the traditional technique using direct estimation of the local
error $\Upsilon^{n-\theta}$,
the stability estimate in \cref{thm: stability} suggests
that one can consider the global consistency error $\mathcal{E}_{\mathrm{glob}}^n$,
defined in \eqref{eq: globalConsistencyError-definition},
accumulated from $t=t_{1-\theta}$ to $t=t_{n-\theta}$ with the complementary discrete
kernel $P_{n-j}^{(n)}$. To exploit this convolution structure,
we will control $\Upsilon^{n-\theta}$
by an \textbf{ECS} bound in terms of the discrete
kernels $A_{n-k}^{(n)}$ defined in \eqref{eq: weights},
and the following quantities
\begin{align}
\Gloc^k&\defeq \frac32\int_{t_{k-1}}^{t_{k-1/2}}(s-t_{k-1})^2|v'''(s)|\zd{s}
	+\frac{3\tau_k}{2}\int_{t_{k-1/2}}^{t_k}(t_k-s)|v'''(s)|\zd{s},\label{eq: Gn loc}\\
\Ghis^k&\defeq\frac52\int_{t_{k-1}}^{t_k}(s-t_{k-1})^2|v'''(s)|\zd{s}
	+\frac52\int_{t_k}^{t_{k+1}}(t_{k+1}-s)^2|v'''(s)|\zd{s},\label{eq: Gn his}
\end{align}
assuming in what follows that $v$ is such that these integrals exist and are
finite.

\subsection{Global consistency error}

\begin{lemma}\label{lem: Upsilon loc}
For any function $v\in C^3((0,T])$, the local consistency
error~$\Upsilon_n^{n-\theta}$ in \eqref{eq: Upsilon n n} satisfies
\[\absb{\Upsilon_n^{n-\theta}}\leq a_0^{(n)}G_{\mathrm{loc}}^n\leq
A_0^{(n)}G_{\mathrm{loc}}^n\quad \text{for~$ 1\leq n\leq N$.}\]
\end{lemma}
\begin{proof}
Taylor expansion (with integral remainder) about~$t_{n-1/2}$ shows that
\[
v'(s)=v'(t_{n-1/2})+v''(t_{n-1/2})(s-t_{n-1/2})
	+\int_{t_{n-1/2}}^s(s-y)v'''(y)\zd{y},
\]
and
\begin{align*}
\brab{\Err{1,n}v}'(s)=&\,v''(t_{n-1/2})(s-t_{n-1/2})
	+\int_{t_{n-1/2}}^s(s-y)v'''(y)\zd{y}\\
	&\,-\frac{1}{2\tau_n}\int_{t_{n-1}}^{t_{n-1/2}}(y-t_{n-1})^2v'''(y)\zd{y}
	-\frac{1}{2\tau_n}\int_{t_{n-1/2}}^{t_n}\!\!(t_n-y)^2v'''(y)\zd{y}.
\end{align*}
Inserting these four terms in~\eqref{eq: Upsilon n n} yields the splitting
$\Upsilon^{n-\theta}_n=\sum_{\ell=1}^4\Upsilon^{n-\theta}_{n,\ell}$. After
integrating by parts, we find that
\begin{equation}\label{eq: theta alpha}
\Upsilon^{n-\theta}_{n,1}
    =(\alpha-2\theta)\,\frac{(1-\theta)^{1-\alpha}}{2\Gamma(3-\alpha)}\,
	v''(t_{n-1/2})\tau_n^{2-\alpha},
\end{equation}
which vanishes for~$\theta=\alpha/2$.
For the term $\Upsilon_{n,2}^{n-\theta}$,
we split the integration interval $[t_{n-1},t_{n-\theta}]$ into two parts:
$[t_{n-1},t_{n-1/2}]$ and $[t_{n-1/2},t_{n-\theta}]$.
Since $t_{n-1/2}<t_{n-\theta}<t_{n}$, the second term reads
\begin{align*}
\Upsilon^{n-\theta}_{n,2}=&\int_{t_{n-1}}^{t_{n-\theta}}\varpi_n'(s)\int_{t_{n-1/2}}^s(s-y)v'''(y)\zd{y}\zd{s}\\
=&\int_{t_{n-1}}^{t_{n-1/2}}\!\!\varpi_n'(s)\!\!\int_s^{t_{n-1/2}}\!\!\!(y-s)v'''(y)\zd{y}\zd{s}
	+\!\int_{t_{n-1/2}}^{t_{n-\theta}}\!\!\!\varpi_n'(s)\!\!\int_{t_{n-1/2}}^s\!\!\!\!\!(s-y)v'''(y)\zd{y}\zd{s}.
\end{align*}
Reversing the order of integration, then integrating by parts in the second
term and using $\varpi_n(t_{n-\theta})=0$, we have
\begin{align*}
\Upsilon^{n-\theta}_{n,2}
	&=\int_{t_{n-1}}^{t_{n-1/2}}\!\!v'''(y)\int_{t_{n-1}}^y\!\!(y-s)\varpi_n'(s)\zd{s}\zd{y}
+\int_{t_{n-1/2}}^{t_{n-\theta}}\!\!v'''(y)
    \int_y^{t_{n-\theta}}\!\!(s-y)\varpi_n'(s)\zd{s}\zd{y}\\
	&=\int_{t_{n-1}}^{t_{n-1/2}}v'''(y)\int_{t_{n-1}}^y(y-s)\varpi_n'(s)\zd{s}\zd{y}
-\int_{t_{n-1/2}}^{t_{n-\theta}}v'''(y)\int_y^{t_{n-\theta}}\varpi_n(s)\zd{s}\zd{y}.
\end{align*}
The inner integrals can be estimated as
\begin{align*}
&\biggl|\int_{t_{n-1}}^y(y-s)\varpi_n'(s)\zd{s}\biggr|
	\le \varpi_n'(t_{n-1/2})\frac{(y-t_{n-1})^2}{2}\quad\text{for $t_{n-1}<y<t_{n-1/2}$,}\\
&\biggl|\int_y^{t_{n-\theta}}\varpi_n(s)\zd{s}\biggr|\le\absb{\varpi_n(t_{n-1/2})}(t_{n-\theta}-y)
	\quad\text{for $t_{n-1/2}<y<t_{n-\theta}$.}
\end{align*}
Recalling the definition \eqref{eq: an} of $a^{(n)}_0$, we see that
$\omega_{2-\alpha}(t_{n-\theta}-t_{n-1})=\tau_na^{(n)}_0$ and then
\begin{align*}
\absb{\varpi_n(t_{n-1/2})}&=\omega_{2-\alpha}(t_{n-\theta}-t_{n-1/2})
	\leq\omega_{2-\alpha}(t_{n-\theta}-t_{n-1})=\tau_na^{(n)}_0,\\
\varpi_n'(t_{n-1/2})&=\omega_{1-\alpha}(t_{n-\theta}-t_{n-1/2})
	=\frac{2}{\tau_n}\omega_{2-\alpha}(t_{n-\theta}-t_{n-1/2})
    \leq 2a^{(n)}_0,
\end{align*}
where we used the fact that $t_{n-\theta}-t_{n-1/2}=(1-\alpha)\tau_n/2$.
Hence, it follows that
\[
\absb{\Upsilon^{n-\theta}_{n,2}}\le a^{(n)}_0\int_{t_{n-1}}^{t_{n-1/2}}
    (y-t_{n-1})^2|v'''(y)|\zd{y}+a^{(n)}_0\tau_n
    \int_{t_{n-1/2}}^{t_{n-\theta}}(t_{n-\theta}-y)|v'''(y)|\zd{y},
\]
and finally
\[
\absB{\sum_{\ell=3}^4\Upsilon^{n-\theta}_{n,\ell}}
	\leq\frac{a^{(n)}_0}{2}
	\int_{t_{n-1}}^{t_{n-1/2}}(y-t_{n-1})^2|v'''(y)|\zd{y}+
	\frac{a^{(n)}_0}{2}\int_{t_{n-1/2}}^{t_n}\!\!(t_n-y)^2|v'''(y)|\zd{y}.
\]
Thus the triangle inequality yields $|\Upsilon^{n-\theta}_n|\le a^{(n)}_0\Gloc^n$
where $\Gloc^n$ is defined in \eqref{eq: Gn loc}.
The definition \eqref{eq: weights} implies $a^{(n)}_0\leq A^{(n)}_0$ and completes the proof.
\end{proof}

\begin{remark}\label{rem: theta}
If we were to choose $\theta\neq\alpha/2$, the
term~\eqref{eq: theta alpha} would limit the consistency error
to an order of~$O(\tau_n^{2-\alpha})$, even for smooth solutions.
\end{remark}

To estimate the remaining terms in~\eqref{eq: Upsilon n k}, we present an
interpolation error formula for the quadratic polynomial $\Pi_{2,k}v$
employed in the Alikhanov formula \eqref{eq: D tau},
but leave the proof to Appendix~\ref{appendix:quadraticInterpolationError}.
This formula is crucial for verifying the \textbf{ECS} hypothesis for the local
consistency error $\Upsilon^{n-\theta}$.

\begin{lemma}\label{lem:quadraticInterpolationError}
If $v\in C^3([t_{k-1},t_{k+1}])$~and $q\in C^2([t_{k-1},t_k])$, then
\begin{align*}
\int_{t_{k-1}}^{t_k}q'(t)\brab{\widetilde{\Pi_{2,k}}v}'(t)\zd{t} =&\,
\int_{t_{k}}^{t_{k+1}}(t_{k+1}-s)^2v'''(s)\zd{s}
\int_{t_{k-1}}^{t_k}\frac{\brab{\widetilde{\Pi_{1,k}}q}(t)\zd t}{(\tau_{k+1}+\tau_{k})\tau_{k+1}}\nonumber\\
&\,-\int_{t_{k-1}}^{t_{k}}(s-t_{k-1})^2v'''(s)\zd{s}
\int_{t_{k-1}}^{t_k}\frac{\brab{\widetilde{\Pi_{1,k}}q}(t)\zd t}{(\tau_{k+1}+\tau_{k})\tau_{k}}\nonumber\\
&\,+\int_{t_{k-1}}^{t_{k}}v'''(s)\zd{s}\int_{t_{k-1}}^{s}\brab{\widetilde{\Pi_{1,k}}q}(t)\zd t,\quad 1\leq k\leq n-1.
\end{align*}
\end{lemma}


\begin{theorem}\label{thm: Upsilon global}
Assume that the mesh condition \Mss{1} holds and
$v\in C^3((0,T])$. For the nonuniform Alikhanov formula~\eqref{eq: D tau} with
the discrete kernels~\eqref{eq: weights}, an \textbf{\upshape ECS} dominates
the local consistency error $\Upsilon^{n-\theta}$ in~\eqref{eq: Upsilon}, that is,
\begin{align*}
\absb{\Upsilon^{n-\theta}}\leq A_{0}^{(n)}\Gloc^n+
\sum_{k=1}^{n-1}\brab{A_{n-k-1}^{(n)}-A_{n-k}^{(n)}}\Ghis^k\quad \text{for $1\leq n\leq N$},
\end{align*}
and consequently the global consistency error satisfies
\begin{align*}
\mathcal{E}_{\mathrm{glob}}^n\leq
\sum_{k=1}^nP_{n-k}^{(n)}A_0^{(k)}\Gloc^k+
\sum_{k=1}^{n-1}P_{n-k}^{(n)}A_{0}^{(k)}\Ghis^k\quad \text{for $1\leq n\leq N$},
\end{align*}
where $\Gloc^k$ and $\Ghis^k$ are defined by \eqref{eq: Gn loc}~and
\eqref{eq: Gn his}, respectively.
\end{theorem}
\begin{proof}
According to the arguments in \cref{subsec: ECS}, it suffices to
verify the first inequality (the \textbf{ECS} bound). The definition 
\eqref{eq:Upsilon n k} of $\Upsilon^{n-\theta}_k$ and
\cref{lem:quadraticInterpolationError} (taking $q\defeq \varpi_n$) yield
\begin{multline}\label{eq: Upsilon his-proof1}
\Upsilon^{n-\theta}_k
=\frac{b^{(n)}_{n-k}}{2}\int_{t_{k-1}}^{t_k}(s-t_{k-1})^2v'''(s)\zd{s}
-\frac{\rho_kb^{(n)}_{n-k}}{2}\int_{t_k}^{t_{k+1}}(t_{k+1}-s)^2v'''(s)\zd{s}\\
	+\int_{t_{k-1}}^{t_k}v'''(s)\int_{t_{k-1}}^s
		\bigl(\Err{1,k}\varpi_n\bigr)(t)\zd{t}\zd{s},\quad 1\leq k\leq n-1,
\end{multline}
where the alternative definition \eqref{eq: bn-anlter} of $b^{(n)}_{n-k}$ has been used.
Recall the error formula of linear interpolation \cite[Lemma 3.1]{LiaoLiZhang:2018},
\begin{align*}
\bigl(\Err{1,k}\varpi_n\bigr)(t)=\int_{t_{k-1}}^{t_k}\chi_k(t,y)\varpi_n''(y)\zd{y},\quad t_{k-1}<t<t_k,\;1\leq k\leq n-1,
\end{align*}
where the Peano kernel
$\chi_k(t,y)=\max\{t-y,0\}-(t-t_{k-1})(t_{k}-y)/\tau_k$ satisfies
\[
-\frac{t-t_{k-1}}{\tau_k}(t_{k}-y)\leq \chi_k(t,y)<0
    \quad\text{for any $t,y\in(t_{k-1},t_k)$.}
\]
The inner integral in the last term of \eqref{eq: Upsilon his-proof1} can be bounded by
\[
\biggl|\int_{t_{k-1}}^s\bigl(\Err{1,k}\varpi_n\bigr)(t)\zd{t}\biggr|
\le\frac{1}{2}(s-t_{k-1})^2\int_{t_{k-1}}^{t_k}\frac{t_{k}-s}{\tau_k}\,
    \varpi_n''(s)\zd{s},\quad t_{k-1}<s<t_k.
\]
By the definition \eqref{eq: Gn his} of $\Ghis^n$ and the triangle inequality, we obtain from \eqref{eq: Upsilon his-proof1} that
\begin{align*}
\absb{\Upsilon^{n-\theta}_k}\le\frac15\braB{(1+\rho_k) b^{(n)}_{n-k}
+\int_{t_{k-1}}^{t_k}\!\!\!\frac{t_{k}-s}{\tau_k}\varpi_n''(s)\zd{s}}\Ghis^k
\leq\brab{A_{n-k-1}^{(n)}-A_{n-k}^{(n)}}\Ghis^k,
\end{align*}
where \cref{thm: A1 A2} (II) was used in the second inequality.
Then the definition \eqref{eq: Upsilon} and \cref{lem: Upsilon loc}
yield the first inequality immediately. The proof is completed.
\end{proof}
\begin{remark}\label{rem: local estimate}
Traditionally, the global approximation error would be estimated by using the truncation error $\Upsilon^{n-\theta}$
directly. Once an upper bound of $\absb{\Upsilon^{n-\theta}}$ is available,
the inequality \eqref{eq: P bound} with $m=0$ will give the global approximate error
\begin{align*}
\mathcal{E}_{\mathrm{glob}}^n\leq\sum_{j=1}^nP_{n-j}^{(n)}\omega_{1-\alpha}(t_j)\max_{1\leq l\leq n}\frac{\abs{\Upsilon^{l-\theta}}}{\omega_{1-\alpha}(t_l)}\leq\pi_A\Gamma(1-\alpha)\max_{1\leq l\leq n}t_l^{\alpha}\absb{\Upsilon^{l-\theta}}\,.
\end{align*}
Nonetheless, the local and global consistency errors described in Theorem \ref{thm: Upsilon global}
present a new understanding of the error contributions generated by the two
different polynomial approximations, respectively, in the local cell $[t_{n-1},t_{n-\theta}]$
and the historical interval $[0,t_{n-1}]$
of the fractional Caputo derivative.

Originally, our \textbf{ECS} bound for~$\Upsilon^{n-\theta}$ is constructed to
preserve the convolution structure of the Caputo fractional derivative as much
as possible.  A direct estimate of the global consistency
error~\eqref{eq: globalConsistencyError-definition} would lead to the double
sum~$\sum_{k=1}^n\sum_{j=1}^k|\Upsilon^{k-\theta}_j|$, whereas the \textbf{ECS}
bound leads to a single
sum~$\sum_{k=1}^nP_{n-k}^{(n)}A_0^{(k)}\brab{\Gloc^k+\Ghis^k}$.  This
simplification assists for proving sharp error bounds even with quite general
nonuniform meshes.
Nonetheless, an explicit bound for the
complementary discrete kernel~$P^{(n)}_{n-j}$ remains an open problem
until now,  and we will make full use of the
identity~\eqref{eq: complementaryKernel-equality} and the
upper bound~\eqref{eq: P bound} in the subsequent analysis.
\end{remark}

\begin{lemma}\label{lem: Upsilon global singular}
Assume that $v\in C^3((0,T])$,
and there exists a positive constant $C_v$ such that $\absb{v'''(t)}\leq C_v(1+t^{\sigma-3})$ for $0<t\leq T$,
where $\sigma\in(0,1)\cup(1,2)$ is a regularity parameter.
If the mesh condition  \Mss{1} holds, for $1\leq n\leq N$, then the global
consistency error satisfies
\begin{align*}
\mathcal{E}_{\mathrm{glob}}^n
\le C_v\braB{\tau_1^{\sigma}/\sigma+t_1^{\sigma-3}\tau_2^3
+\frac1{1-\alpha}\max_{2\leq k\leq n}t_{k}^{\alpha}t_{k-1}^{\sigma-3}\tau_k^{3}/\tau_{k-1}^{\alpha}}.
\end{align*}
\end{lemma}
\begin{proof}
The bounds on the discrete kernel~$A_{n-k}^{(n)}$ in \cref{thm: A1 A2} (I) yield
the inequalities
$$A_{0}^{(k)}\leq \frac{24}{11}\omega_{2-\alpha}(\tau_k)/\tau_k,\quad
A_{k-2}^{(k)}\geq\frac4{11}\omega_{1-\alpha}(t_k-t_1),$$
and
\begin{align*}
\frac{A_{0}^{(k)}}{A_{k-2}^{(k)}}<\frac{6\,\omega_{2-\alpha}(\tau_k)}{\tau_k\,\omega_{1-\alpha}(t_k-t_1)}
\leq\frac{6}{1-\alpha}\frac{(t_k-t_1)^{\alpha}}{\tau_k^{\alpha}},\quad 2\leq k\leq n\le N.
\end{align*}
Furthermore, the identity~\eqref{eq: complementaryKernel-equality} for the complementary discrete
kernel~$P^{(n)}_{n-j}$ gives
\[
P_{n-1}^{(n)}A_{0}^{(1)}\leq1\quad \mbox{and}\quad
\sum_{k=2}^{n-1}P_{n-k}^{(n)}A_{k-2}^{(k)}
    \le\sum_{k=2}^{n}P_{n-k}^{(n)}A_{k-2}^{(k)}=1.
\]
Applying the definition \eqref{eq: Gn loc} with the regularity assumption,
it is not difficult to get
$$\Gloc^1\le C_v\tau_1^{\sigma}/\sigma\quad \text{and}\quad
\Gloc^k\leq C_v t_{k-1}^{\sigma-3}\tau_k^{3}\quad\text{ for $2\leq k\leq N.$}$$
Similarly, by using the formula~\eqref{eq: Gn his}, one gets
$$\Ghis^1\leq
C_v(\tau_1^{\sigma}/\sigma+t_1^{\sigma-3}\tau_2^3)\;\;\text{and}\;\;
\Ghis^k\leq C_v \bra{t_{k-1}^{\sigma-3}\tau_k^{3}
+t_{k}^{\sigma-3}\tau_{k+1}^{3}}\;\;\text{for $2\leq k\leq N-1$.}$$
Then it follows from \cref{thm: Upsilon global}  that
\[
\mathcal{E}_{\mathrm{glob}}^n
\leq P_{n-1}^{(n)}A_{0}^{(1)}\brab{\Gloc^1+\Ghis^1}
    +\sum_{k=2}^{n}P_{n-k}^{(n)}A_{0}^{(k)}\Gloc^k
    +\sum_{k=2}^{n-1}P_{n-k}^{(n)}A_{0}^{(k)}\Ghis^k.
\]
The first term on the right is bounded
by~$C_v(\tau_1^\sigma/\sigma+t_1^{\sigma-3}\tau_2^3)$, and the remaining
terms can be bounded by
\begin{align*}
&\frac{6}{1-\alpha}\biggl(\sum_{k=2}^{n}P_{n-k}^{(n)}A_{k-2}^{(k)}
    t_k^{\alpha}\tau_k^{-\alpha}\Gloc^k
+\sum_{k=2}^{n-1}P_{n-k}^{(n)}A_{k-2}^{(k)}
    t_k^{\alpha}\tau_k^{-\alpha}\Ghis^k\biggr)\\
&\leq\frac{C_v}{1-\alpha}\max_{2\leq k \leq n}
    t_k^{\alpha} t_{k-1}^{\sigma-3}\tau_k^{3-\alpha}
    +\frac{C_v}{1-\alpha}\max_{2\leq k\leq n-1}
    \Bigl(t_k^{\alpha}t_{k-1}^{\sigma-3}\tau_k^{3-\alpha}
    + t_k^{\alpha+\sigma-3}\tau_{k+1}^{3}\tau_k^{-\alpha}\Bigr) \\
    &\leq\frac{C_v}{1-\alpha}\max_{2\leq k \leq n}
    t_k^\alpha t_{k-1}^{\sigma-3}
\tau_k^{3}/\tau_{k-1}^{\alpha}(1+\rho_{k-1}^{\alpha}),
\end{align*}
implying the claimed estimate.
\end{proof}

\begin{remark}\label{rem: firstlevel}
The proof of Lemma \ref{lem: Upsilon global singular} and the \textbf{ECS} bound in Theorem \ref{thm: Upsilon global} give
\[
\absb{\Upsilon^{1-\theta}}\leq A_{0}^{(1)}\Gloc^1\leq
C_v\tau_1^{\sigma-\alpha}/\sigma,
\]
implying that $\Upsilon^{1-\theta}=O(1)$ when $\sigma=\alpha$, and
if $0<\sigma<\alpha$ then the situation becomes worse. The global consistency
analysis seems therefore to be also a superconvergence analysis.
\end{remark}

Now we describe the contribution to the global truncation error from the time
weighted terms. The next lemma suggests that the temporal error introduced by
the time weighted approach is smaller than that generated by the Alikhanov
approximation of the Caputo derivative.

\begin{lemma}\label{lem:time-Weighted global}
Assume that $v\in C^2((0,T])$,
and there exists a positive constant $C_v$ such that $\absb{v''(t)}\leq C_v(1+t^{\sigma-2})$ for $0<t\leq T$,
where $\sigma\in(0,1)\cup(1,2)$ is a regularity parameter.
Denote the local truncation error of $v^{n-\theta}$ by
\begin{align*}
\mathcal{R}^{n-\theta}=v(t_{n-\theta})-v^{n-\theta}\quad \text{for $1\leq n\leq N$.}
\end{align*}
If the mesh condition \Mss{1} holds, then the global consistency error satisfies
\begin{align*}
\sum_{j=1}^nP_{n-j}^{(n)}\absb{\mathcal{R}^{j-\theta}}
\leq&\,C_v\braB{\tau_1^{\sigma+\alpha}/\sigma
+t_n^{\alpha}\max_{2\leq k\leq n}t_{k-1}^{\sigma-2}\tau_k^{2}}\quad \text{for $1\leq n\leq N$.}
\end{align*}
\end{lemma}
\begin{proof}
The following integral representation of $\mathcal{R}^{j-\theta}$ can be easily
verified, for example using the Taylor formula with integral
remainder~\cite[Lemma~2.5]{LiaoZhaoTeng:2016},
\begin{align*}
\mathcal{R}^{j-\theta}=-\theta\int_{t_{j-1}}^{t_{j-\theta}}(s-t_{j-1})v''(s)\zd{s}
-(1-\theta)\int^{t_{j}}_{t_{j-\theta}}(t_{j}-s)v''(s)\zd{s}\,,\quad 1\leq j\leq N.
\end{align*}
Under the regularity assumption,
one has
\[
\absb{\mathcal{R}^{1-\theta}}\leq C_v\,\frac{\tau_1^\sigma}{\sigma}
\quad\text{and}\quad
\absb{\mathcal{R}^{j-\theta}}\leq C_vt_{j-1}^{\sigma-2}\tau_j^{2},
\quad 2\leq j\leq N.
\]
Note that \cref{thm: A1 A2} (I) implies $A_{0}^{(1)}\geq \frac4{11}
\omega_{2-\alpha}(\tau_1)/\tau_1$,
and then the identity~\eqref{eq: complementaryKernel-equality} shows that
$$P_{n-1}^{(n)}\leq1/A_{0}^{(1)}\leq 3\Gamma(2-\alpha)\tau_1^{\alpha}.$$
Therefore we obtain
\begin{align*}
\sum_{j=1}^nP_{n-j}^{(n)}\absb{\mathcal{R}^{j-\theta}}=&\,P_{n-1}^{(n)}\absb{\mathcal{R}^{1-\theta}}
+\sum_{j=2}^nP_{n-j}^{(n)}\absb{\mathcal{R}^{j-\theta}}\\
\leq&\,3\Gamma(2-\alpha)\tau_1^{\alpha}\absb{\mathcal{R}^{1-\theta}}
+\max_{2\leq k\leq n}\absb{\mathcal{R}^{k-\theta}}\sum_{j=1}^nP_{n-j}^{(n)}\\
\leq&\,C_v\braB{\tau_1^{\sigma+\alpha}/\sigma
+t_n^{\alpha}\max_{2\leq k\leq n}t_{k-1}^{\sigma-2}\tau_k^{2}},
\qquad 1\leq n\leq N,
\end{align*}
where the estimate \eqref{eq: P bound} with $\pi_A=11/4$ has been used in the last inequality.
\end{proof}

\subsection{Convergence}

We now establish the convergence of the numerical solution under the regularity
conditions~\eqref{eq: sigma} and the assumptions \Mss{1}--\Mss{2}.
To deal with the spatial error, we introduce the Ritz
projector~$R_h:H^1_0(\Omega)\to X_h$, defined by
\[
\iprod{\nabla R_hv,\nabla\chi}=\iprod{\nabla v,\nabla\chi}
	\quad\text{for $v\in H^1_0(\Omega)$ and $\chi\in X_h$.}
\]
\begin{theorem}\label{thm: convergence}
Suppose that the solution $u$ of \eqref{eq: IBVP} has the regularity
property \eqref{eq: sigma} for the parameter~$\sigma\in(0,1)\cup(1,2)$, and
consider the time-stepping method \eqref{eq: discrete IBVP} using the nonuniform
Alikhanov formula \eqref{eq: D tau} with the discrete
kernels~\eqref{eq: weights}. If \Mss{1} holds with the
maximum step size $\taumax\le1/\sqrt[\alpha]{11\Gamma(2-\alpha)\kappa_{+}}$, then
the discrete solution $u_h^n$ is convergent with respect to the $L_2$-norm,
\begin{multline*}
\mynormb{u(t_n)-u^n_h}
	\le C_uE_\alpha(20\kappa_{+}t_n^\alpha) \bigg(
\frac{\tau_1^\sigma}{\sigma}
	+\frac{1}{1-\alpha}\max_{2\le k\le n}
	t_k^\alpha t_{k-1}^{\sigma-3}\,\frac{\tau_k^{3}}{\tau_{k-1}^{\alpha}}
+t_n^\alpha\max_{2\le k\le n}t_{k-1}^{\sigma-2}\tau_k^2	\\
    +\|u_{0h}-R_hu_0\| +(t_n+t_n^{\alpha}+t_n^\sigma)h^2 \bigg)
\quad \text{for $1\leq n\leq N$}.
\end{multline*}
In particular, if \Mss{2} also holds  {and if we choose $u_{0h}=R_hu_0$}, then
\[
\mynormb{u(t_n)-u^n_h}\le
\frac{C_u}{\sigma(1-\alpha)}\,\tau^{\min\{\gamma\sigma,2\}}+C_uh^2
	\quad\text{for $1\le n\le N$,}
\]
where $C_u$ may depend on $u$ and $T$, but is uniformly bounded with respect to $\alpha$~and $\sigma$.
\end{theorem}
\begin{proof}
Let $e^n_h=u^n_h-R_hu^n\in X_h$ where $u^n=u(t_n)$, so that
\[
\|u^n_h-u^n\|\le\|u^n-R_hu^n\|+\|e^n_h\|.
\]
The usual analysis of the elliptic problem shows that, under the first
regularity assumption in~\eqref{eq: sigma},
\begin{equation}\label{eq: Ritz error}
\|u^n-R_hu^n\|\le C_{\Omega}h^2\|u^n\|_{H^2(\Omega)}\le C_uh^2,
\end{equation}
so it suffices to deal with~$e^n_h$.  We find
\cite[Section~4]{LiaoMcLeanZhang:2019} that
\[
\iprod{(\dfd{\alpha}e_h)^{n-\theta},\chi}
	+\iprod{\nabla e_h^{n-\theta},\nabla\chi}
        =\kappa\iprod{e_h^{n-\theta},\chi}+\iprod{\mathcal{R}^n,\chi},
\]
for all $\chi\in X_h$, where
\begin{equation}\label{eq: Rn def}
\mathcal{R}^n=(\fd{\alpha}u)(t_{n-\theta})-(\dfd{\alpha}R_hu)^{n-\theta}
	-\kappa\bigl(u(t_{n-\theta})-R_hu^{n-\theta}\bigr)
	+\triangle\bigl(u^{n-\theta}-u(t_{n-\theta})\bigr).
\end{equation}
Choosing $\chi=u^{n-\theta}_h$ yields an inequality
of the form~\eqref{eq: energy} with $u^{n-\theta}_h$~and $f(t_{n-\theta})$
replaced by $e^{n-\theta}_h$~and $\mathcal{R}^n$, respectively.
Hence, the argument leading to Theorem~\ref{thm: stability} shows that
\begin{equation}\label{eq: e R}
\|e^n_h\|\le2E_\alpha(20\kappa_{+}t_n^\alpha)\biggl(\|e^0_h\|
	+2\max_{1\le k\le n}\sum_{j=1}^k P^{(k)}_{k-j}
		\|\mathcal{R}^j\|\biggr)\quad\text{for $1\le n\le N$.}
\end{equation}
Write $\mathcal{R}^j=\mathcal{R}^j_1+\mathcal{R}^j_2
+\mathcal{R}^j_3+\mathcal{R}^j_4$, where
\begin{align*}
\mathcal{R}^j_1&=
	(\fd{\alpha}u)(t_{j-\theta})-(\dfd{\alpha}u)^{j-\theta},&
\mathcal{R}^j_2&=
	(\kappa+\triangle)\bigl(u^{j-\theta}-u(t_{j-\theta})\bigr),\\
\mathcal{R}^j_3&=\bigl(\dfd{\alpha}(u-R_hu)\bigr)^{j-\theta},&
\mathcal{R}^j_4&=\kappa(R_hu-u)^{j-\theta}.
\end{align*}
Applying \cref{lem: Upsilon global singular,lem:time-Weighted global},
combined with the regularity assumption~\eqref{eq: sigma},
one obtains
\[
\max_{1\le k\le n}\sum_{j=1}^k P^{(k)}_{k-j}
	\bigl\|\mathcal{R}^j_1+\mathcal{R}^j_2\bigr\|
\le C_u\biggl(\frac{\tau_1^{\sigma}}{\sigma}
    +\frac1{1-\alpha}\max_{2\leq k\leq n}t_k^\alpha
		t_{k-1}^{\sigma-3}\,\frac{\tau_k^{3}}{\tau_{k-1}^{\alpha}}
    +t_n^{\alpha}\max_{2\leq k\leq n}t_{k-1}^{\sigma-2}\tau_k^{2}\biggr).
\]
Since
\[
\|\mathcal{R}^j_3\|=\biggl\|\sum_{\ell=1}^jA^{(j)}_{j-\ell}
	\nabla_\tau(u-R_hu)^\ell\biggr\|
	\le\sum_{\ell=1}^j A^{(j)}_{j-\ell}\int_{t_{\ell-1}}^{t_\ell}
		\|(u-R_hu)'(t)\|\,dt,
\]
the identity~\eqref{eq: complementaryKernel-equality}, the error bound~\eqref{eq: Ritz error}
for the Ritz projection and the regularity assumption~\eqref{eq: sigma} give
\begin{align*}
\max_{1\le k\le n}\sum_{j=1}^k P^{(k)}_{k-j}\|\mathcal{R}^j_3\|
	&\le\max_{1\le k\le n}\sum_{\ell=1}^k
		\biggl(\sum_{j=\ell}^k P^{(k)}_{k-j} A^{(j)}_{j-\ell}\biggr)
		\int_{t_{\ell-1}}^{t_\ell}\|(u-R_hu)'(t)\|\,dt\\
	&\le C_uh^2\int_0^{t_n}\|u'(t)\|_{H^2(\Omega)}\,dt
	\le C_u(t_n+t_n^\sigma)h^2.
\end{align*}
Recalling the upper bound~\eqref{eq: P bound}~and the Ritz projection error~\eqref{eq: Ritz error}, we see that
\[
\max_{1\le k\le n}\sum_{j=1}^k P^{(k)}_{k-j}\|\mathcal{R}^j_4\|
	\le C_{\Omega}h^2\max_{1\le k\le n}\sum_{j=1}^k P^{(k)}_{k-j}
		\|u^{j-\theta}\|_{H^2(\Omega)}
	\le C_u t_n^\alpha h^2,
\]
so the first estimate for~$\|u^n_h-u(t_n)\|$ follows.
If the mesh assumption~$\Mss{2}$ holds, then
$\tau_1\le C_{\gamma}\tau^{\gamma}$ and, with~$\beta:=\min\{2,\gamma\sigma\}$,
\begin{align}\label{eq: D tau order}
t_{k}^{\alpha}t_{k-1}^{\sigma-3}\tau_k^{3}/\tau_{k-1}^{\alpha}&\leq C_{\gamma}t_{k}^{\alpha+\sigma-3}\tau_k^{3-\alpha}
\leq C_\gamma t_k^{\sigma-3+\alpha}\tau_k^{3-\alpha-\beta}
    \bigl(\tau\min\{1,t_k^{1-1/\gamma}\}\bigr)^\beta\\
&\leq C_\gamma t_k^{\sigma-\beta/\gamma}(\tau_k/t_k)^{3-\alpha-\beta}\tau^\beta\nonumber\\
&\leq C_\gamma t_k^{\max\{0,\sigma-(3-\alpha)/\gamma\}}\tau^{\beta},
    \quad 2\leq k\leq n. \nonumber
\end{align}
In addition, we have
\begin{multline}\label{eq: offset order}
t_{k-1}^{\sigma-2}\tau_k^2
    \leq C_\gamma t_k^{\sigma-2}\tau_k^{2-\beta}
        \bigl(\tau\min\{1,t_k^{1-1/\gamma}\}\bigr)^\beta\\
\le C_\gamma t_k^{\sigma-\beta/\gamma}(\tau_k/t_k)^{2-\beta}\tau^\beta
\le C_\gamma t_k^{\max\{0,\sigma-2/\gamma\}}\tau^\beta,\quad 2\leq k\le n,
\end{multline}
so the claimed second result follows immediately by noting that $t_n\leq T$.
\end{proof}

\begin{remark}
Replacing $f(t_{n-\theta})$ with~$f^{n-\theta}$
in~\eqref{eq: discrete IBVP} would introduce an additional
term~$f^{n-\theta}-f(t_{n-\theta})$ in the
definition~\eqref{eq: Rn def} of~$\mathcal{R}^n$, but would not affect the
final error bound, assuming $f$ has the regularity properties needed to apply
\cref{lem:time-Weighted global}.
Also, instead of $u_{0h}=R_hu_0$ we could choose the
interpolant or the $L_2$-projection of~$u_0$ and still maintain
second-order accuracy in space.
\end{remark}

\begin{remark}
By an argument similar to that in \eqref{eq: D tau order}, it is not difficult
to show
$$t_{k}^{\alpha}t_{k-1}^{\sigma-3}\tau_k^{3}/\tau_{k-1}^{\alpha}
\le C_{\gamma}t_k^{\sigma-(3-\alpha)/\gamma}\tau^{3-\alpha},$$
which means that the Alikhanov formula $(\dfd{\alpha}v)^{n-\theta}$ approximates
$(\fd{\alpha}u)(t_{n-\theta})$ to order~$O(\tau^{3-\alpha})$ if
$\gamma\ge(3-\alpha)/\sigma$.  However, the term~\eqref{eq: offset order}
arising from the difference~$u(t_{n-\theta})-u^{n-\theta}$
in~\eqref{eq: Rn def} would still limit the convergence rate for the
overall scheme to order~$O(\tau^{2})$.
\end{remark}

\section{Proof of \cref{thm: A1 A2} (discrete convolution kernels)}
\label{sec: monotone}

Our aim is to prove the boundedness and monotonicity of the convolution
kernels~$A_{n-k}^{(n)}$. Since the coefficients $a_{n-k}^{(n)}$,
$b_{n-k}^{(n)}$~and $A_{n-k}^{(n)}$ in \eqref{eq: an}, \eqref{eq: bn}~and
\eqref{eq: weights} are defined on nonuniform meshes, it is a technically
challenging task and some new techniques will be necessary.

Note that, some techniques using Taylor expansion or function monotonicity
have been applied in \cite{Alikhanov2015,GaoSunZhang:2014,LvXu:2016} to investigate
the discrete convolution kernels in high-order numerical Caputo formulas. These techniques
can not be directly applied here although they would be well suited for the uniform case with $\tau_k=\tau$.

We start from the alternative definition \eqref{eq: bn-anlternative} of $b_{n-k}^{(n)}$.
Compared with the original definition \eqref{eq: bn},
the new formula \eqref{eq: bn-anlternative} has a nonnegative integrand
because $\Gamma(-\alpha)<0$ and $\omega_{-\alpha}(t_{n-\theta}-s)<0$,
motivating us to consider the integrals
\begin{align*}
\int_{t_{k-1}}^{t_k}\varpi_n''(t)\zd{t},\quad
\int_{t_{k-1}}^{t_k}\frac{t_{k}-t}{\tau_k}\varpi_n''(t)\zd{t}
\quad\text{and}\quad
\int_{t_{k-1}}^{t_k}\frac{t-t_{k-1}}{\tau_k}\varpi_n''(t)\zd{t}\,.
\end{align*}
Actually, they are very close to the values of $b_{n-k}^{(n)}$ and $a_{n-k-1}^{(n)}-a_{n-k}^{(n)}$
if $\rho=O(1)$. In studying theoretical properties of the discrete kernels
$A^{(n)}_{n-k}$, these integrals will play a bridging role in establishing some
useful links between the underlaying discrete coefficients $a_{n-k}^{(n)}$,
$b_{n-k}^{(n)}$~and $A_{n-k}^{(n)}$; see
\cref{lem: bn qn,lem: bn an,lem: ank diff}.

\subsection{Proof of \cref{thm: A1 A2} (I)}
\begin{lemma}\label{lem: an}
The discrete coefficients~$a^{(n)}_{n-k}$ defined in \eqref{eq: an} satisfy
\begin{itemize}
\item[(i)]
$\displaystyle a^{(n)}_{n-k}>\omega_{1-\alpha}(t_{n-\theta}-t_{k-1})>a^{(n)}_{n-k+1}$ for $1\le k\le n;$
\item[(ii)]
  $ a^{(n)}_0>\frac3{4}\int_{t_{n-1}}^{t_n}\frac{\omega_{1-\alpha}(t_n-s)}{\tau_n}\zd{s}$ and
$a^{(n)}_{n-k}>\int_{t_{k-1}}^{t_k}\frac{\omega_{1-\alpha}(t_n-s)}{\tau_k}\zd{s}$ for $1\le k\le n-1$.
\end{itemize}
\end{lemma}
\begin{proof}(i) If $k=n$, one has
$$a^{(n)}_{0}=\frac{1-\theta}{1-\alpha}\omega_{1-\alpha}(t_{n-\theta}-t_{n-1})>\omega_{1-\alpha}(t_{n-\theta}-t_{n-1}).$$
For $1\leq k<n$, the claimed inequalities follow directly from the integral mean value theorem and
the fact that $\varpi_n'(s)=\omega_{1-\alpha}(t_{n-\theta}-s)$ is a strictly increasing function.

(ii) Also, the lower bounds of $a_{n-k}^{(n)}$ for $1\le k< n$ follow from
the definition \eqref{eq: an} immediately. For the remaining coefficient
$a_{0}^{(n)}$, since $e^{x}>1+x$ for all real~$x$ and since $\ln(1-x/2)>-x$
for~$0<x<1$, we find that
\[
(1-\theta)^{1-\alpha}=e^{(1-\alpha)\ln(1-\alpha/2)}
>1+(1-\alpha)\ln(1-\alpha/2)>1-\alpha(1-\alpha)\geq3/4,
\]
and then
$a_{0}^{(n)}=(1-\theta)^{1-\alpha}\omega_{2-\alpha}(\tau_n)/\tau_n
>\frac3{4\tau_n}\int_{t_{n-1}}^{t_n}\omega_{1-\alpha}(t_n-s)\zd{s}.$
\end{proof}

\begin{lemma}\label{lem: bn qn}
The discrete coefficients $b^{(n)}_{n-k}$ defined in~\eqref{eq: bn} satisfy
\[
0<b^{(n)}_{n-k}\le\frac{\rho_k}{4(1+\rho_k)}\int_{t_{k-1}}^{t_k}\varpi_n''(s)\zd{s},
    \quad1\le k\le n-1.
\]
\end{lemma}
\begin{proof}
Since $0<(s-t_{k-1})(t_k-s)<\tau_k^2/4$ for~$t_{k-1}<s<t_k$, the alternative definition
\eqref{eq: bn-anlternative} of $b^{(n)}_{n-k}$ yields the result and completes the proof.
\end{proof}

As an application of~\cref{lem: bn qn}, the next lemma builds up a link between
$a_{n-k}^{(n)}$~and $b_{n-k}^{(n)}$.
For a uniform mesh~$t_n=n\tau$, this lemma gives
\[
0<b_{n-k}^{(n)}<\frac{\theta }{4(n-\theta-k)}a^{(n)}_{n-k}
\quad \text{for $1\le k\le n-1.$}
\]
By comparison, Alikhanov~\cite[Lemma~3 and
Corollary~2]{Alikhanov2015} gives
$b_{n-k}^{(n)}<\frac{\theta}{2(1-\theta)}a^{(n)}_{n-k}$.
Obviously, the new bound is much sharper.

\begin{lemma}\label{lem: bn an}
The positive coefficients $a^{(n)}_{n-k}$, $b^{(n)}_{n-k}$ defined in \eqref{eq: an} and \eqref{eq: bn} satisfy
\[
b^{(n)}_{n-k}<\frac{\theta\tau_k}{2(t_{n-\theta}-t_k)}\,\frac{\rho_k}{1+\rho_k}a^{(n)}_{n-k},\quad 1\le k\le n-1.
\]
\end{lemma}
\begin{proof}For fixed $n$ and $1\leq k\leq n-1$, consider an auxiliary function
\[
\varphi_k(z)\defeq\int_{t_{k-1}}^{t_{k-1}+z}\varpi_n''(s)\zd{s}
    -\frac{2\theta}{t_{n-\theta}-t_k}\int_{t_{k-1}}^{t_{k-1}+z}\varpi_n'(s)\zd{s},\quad 0<z<\tau_k.
\]
Since $\varpi_n''(t)=\alpha \varpi_n'(t)/(t_{n-\theta}-t)$, the first derivative
\begin{align*}
\varphi_k'(z)&=\varpi_n'(t_{k-1}+z)\biggl(
    \frac{\alpha}{t_{n-\theta}-t_{k-1}-z}
        -\frac{2\theta}{t_{n-\theta}-t_{k}}\biggr)\\
    &<\varpi_n'(t_{k-1}+z)\frac{\alpha-2\theta}{t_{n-\theta}-t_k}=0,
        \quad0<z<\tau_k,\;1\leq k\leq n-1.
\end{align*}
Hence the definition \eqref{eq: an} of $a_{n-k}^{(n)}$ yields
\[
\int_{t_{k-1}}^{t_{k}}\varpi_n''(s)\zd{s}-\frac{2\theta\tau_k}{t_{n-\theta}-t_k}a^{(n)}_{n-k}
    =\varphi_k(\tau_k)<\varphi_k(0)=0,\quad1\leq k\leq n-1.
\]
\cref{lem: bn qn} gives the claimed inequality and completes the proof.
\end{proof}

Now we verify \cref{thm: A1 A2} (I) by using \cref{lem: an,lem: bn an}.
\begin{proof}[Proof of \cref{thm: A1 A2} (I)]
Under the assumption \Mss{1}, one has
$\theta<1-\theta$, $\rho_k\leq7/4$  and
$t_{n-\theta}-t_k\ge(1-\theta)\tau_{k+1}$ for $1\leq k\leq n-1$.
Thus, by using \cref{lem: bn an}, one has
\[
b^{(n)}_{n-k}
	<\frac{\theta\tau_k}{2(1-\theta)\tau_{k+1}}\, \frac{\rho_k }{1+\rho_k}a^{(n)}_{n-k}
    \le\frac{7\rho_k}{8(1+\rho_k)}a^{(n)}_{n-k}
\le\frac{7}{11}a^{(n)}_{n-k},\quad 1\leq k\leq n-1,
\]
since the function $t/(1+t)$ is increasing for any $t>0$.
By \cref{lem: an} (i), $a_1^{(n)}<a_0^{(n)}$, then the definition \eqref{eq: weights} yields
$$A^{(n)}_0=a^{(n)}_0+\rho_{n-1}b^{(n)}_1\leq a^{(n)}_0+\tfrac{49}{44}a^{(n)}_1\leq\frac{24}{11}a_{0}^{(n)}.$$
So the definition \eqref{eq: an}
of $a_0^{(n)}$ gives the upper bound
$$A^{(n)}_0\leq\frac{24}{11\tau_n}(1-\theta)^{1-\alpha}\int_{t_{n-1}}^{t_{n}}\omega_{1-\alpha}(t_n-s)\zd s \leq\frac{24}{11\tau_n}\int_{t_{n-1}}^{t_{n}}\omega_{1-\alpha}(t_n-s)\zd s.$$
The lower bounds of $A^{(n)}_{n-k}$ for~$1\le k\le n-1$ follow from
\cref{lem: an} (ii) because
the definition \eqref{eq: weights} implies that
$$A^{(n)}_{n-k}\ge a^{(n)}_{n-k}-b^{(n)}_{n-k}\ge \frac{4}{11}a^{(n)}_{n-k}.$$
The proof of \cref{thm: A1 A2} (I) is complete.
\end{proof}

\subsection{Proof of \cref{thm: A1 A2} (II)--(III)}
For simplicity of presentation, in this subsection we let
\begin{equation}\label{ank diff bound-notations}
I_{n-k}^{(n)}:=\int_{t_{k-1}}^{t_{k}}\frac{t_{k}-t}{\tau_{k}}\varpi_n''(t)\zd{t}
\quad\mbox{and}\quad
J_{n-k}^{(n)}:=\int_{t_{k-1}}^{t_{k}}\frac{t-t_{k-1}}{\tau_{k}}
    \varpi_n''(t)\zd{ t}
\end{equation}
for $1\leq k\leq n-1$.
\begin{lemma}\label{lem: bn qn2}
For $1\le k\le n-2$, the positive coefficients $b^{(n)}_{n-k}$ in~\eqref{eq: bn} satisfy
\begin{itemize}
\item[(i)]
$ I_{n-k}^{(n)}\geq \frac{1+\rho_k}{\rho_k}b^{(n)}_{n-k}$;
\quad(ii)\;\;
$ J_{n-k}^{(n)}\geq \frac{2(1+\rho_k)}{\rho_k}b^{(n)}_{n-k}$;
\quad(iii)\;\;
$ J_{n-k}^{(n)}\geq I_{n-k}^{(n)}$.
\end{itemize}
\end{lemma}
\begin{proof}The alternative definition
\eqref{eq: bn-anlternative} of $b^{(n)}_{n-k}$ gives the result (i) directly since
$0<s-t_{k-1}<\tau_k$ for $s\in(t_{k-1},t_k)$.
Since $\varpi_n'''(t)>0$ for $0<t<t_{n-\theta}$,
we take $q:=\varpi_n'$ in Lemma \ref{lem: integralFormula} to find
\[
\int_{t_{k-1}}^{t_k}\bra{\frac{s-t_{k-1}}{\tau_k}-\frac{1}2}\varpi_n''(s)\zd{s}
=\frac{1}{2\tau_k}\int_{t_{k-1}}^{t_k}(s-t_{k-1})(t_k-s)\varpi_n'''(s)\zd{s}>0,
\]
and then $J_{n-k}^{(n)}>\frac12\int_{t_{k-1}}^{t_k}\varpi_n''(s)\zd{s}$ for  $1\le k\le n-1.$
So the inequality (ii) follows immediately from \cref{lem: bn qn}. Moreover,
$2J^{n}_{n-k}>\int_{t_{k-1}}^{t_k}\varpi_n''(s)\zd{s}
=I_{n-k}^{(n)}+J_{n-k}^{(n)}$ so the claimed result~(iii) follows directly.
\end{proof}
\begin{lemma}\label{lem: qn qn1}
For any fixed $n$ $(3\leq n\leq N)$ and $1\le k\le n-2$, it holds that
\begin{itemize}
\item[(i)]
$I_{n-k-1}^{(n)}\geq
\frac1{\rho_k}I_{n-k}^{(n)}$; \quad(ii)\;\; $J_{n-k-1}^{(n)}\geq
\frac1{\rho_{k}}J_{n-k}^{(n)}$.
\end{itemize}
\end{lemma}
\begin{proof}
For fixed~$n\geq2$, introduce an auxiliary function with respect to $z\in[0,1]$,
\[
\psi_k(z)\defeq\frac{1}{\tau_k}\int_{t_{k-1}}^{t_{k-1}+z\tau_k}
	\bra{t_{k-1}+z\tau_k-s}\varpi_n''(s)\zd{s},\quad 1\leq k\leq n-1,
\]
with the first and second derivatives
\begin{align*}
\psi_k'(z)=\int_{t_{k-1}}^{t_{k-1}+z\tau_k}\varpi_n''(s)\zd{s},\quad
\psi_k''(z)=\tau_k\varpi_n''(t_{k-1}+z\tau_k),\quad 1\le k\le n-1.
\end{align*}
Note that $\psi_k(0)=\psi_k'(0)=0$ for $1\le k\le n-1$,
and $\psi_{k+1}(0)=\psi_{k+1}'(0)=0$ for $0\le k\le n-2$.
Thanks to the Cauchy differential
mean-value theorem, there exist $z_{1k}$, $z_{2k}\in(0,1)$ such that
\begin{align*}
\frac{I^{(n)}_{n-k-1}}{I^{(n)}_{n-k}}&=\frac{\psi_{k+1}(1)}{\psi_k(1)}
	=\frac{\psi_{k+1}(1)-\psi_{k+1}(0)}{\psi_k(1)-\psi_k(0)}
	=\frac{\psi_{k+1}'(z_{1k})}{\psi_k'(z_{1k})}
	=\frac{\psi_{k+1}'(z_{1k})-\psi_{k+1}'(0)}{\psi_k'(z_{1k})-\psi_k'(0)}
	\\
	&=\frac{\psi_{k+1}''(z_{2k})}{\psi_k''(z_{2k})}
     =\frac{\tau_{k+1}\varpi_n''(t_{k}+z_{2k}\tau_{k+1})}{\tau_{k}\varpi_n''(t_{k-1}+z_{2k}\tau_k)}\geq \frac1{\rho_{k}},
     \quad 1\leq k\leq n-2,
\end{align*}
because $\varpi_n''(t)>0$ is increasing and
$t_k>t_{k-1}+z_{2k}\tau_k$. The inequality~(i) follows.
We now introduce another auxiliary function for~$z\in[0,1]$,
\[
\phi_k(z)\defeq\frac{1}{\tau_k}\int_{t_{k-1}}^{t_{k-1}+z\tau_k}
	(s-t_{k-1})\varpi_n''(s)\zd{s},\quad 1\leq k\leq n-1,
\]
with the first derivative $\phi_k'(z)=z\tau_k\varpi_n''(t_{k-1}+z\tau_k)$ for $1\le k\le n-1.$
Then a similar argument yields the desired result~(ii) and completes the proof.
\end{proof}

\begin{lemma}\label{lem: ank diff}
The positive coefficients $a^{(n)}_{n-k}$ in \eqref{eq: an} satisfy
 \begin{align*}
a_{n-k-1}^{(n)}-a_{n-k}^{(n)}
=I_{n-k-1}^{(n)}+J_{n-k}^{(n)},\quad 1\le k\le n-2\; (3\leq n\leq N),
 \end{align*}
and for $k=n-1$ $(2\leq n\leq N)$,
 \begin{align*}
a_{0}^{(n)}-a_{1}^{(n)}
=\frac{\theta}{1-2\theta}\varpi_n'(t_{n-1})+J_{1}^{(n)}.
 \end{align*}
\end{lemma}
\begin{proof} For fixed $n$ $(3\leq n\leq N)$, applying the definition \eqref{eq: an},
we exchange the order of integration to find
\begin{align*}
a_{n-k-1}^{(n)}-\varpi_n'(t_k)
=&\,\int_{t_{k}}^{t_{k+1}}\frac{\varpi_n'(s)-\varpi_n'(t_k)}{\tau_{k+1}}\zd s
=\int_{t_{k}}^{t_{k+1}}\!\!\int_{t_{k}}^s\frac{\varpi_n''(t)}{\tau_{k+1}}\zd{t}\zd s=I_{n-k-1}^{(n)}
 \end{align*}
for $0\leq k\leq n-2$, and similarly,
\begin{align*}
a_{n-k}^{(n)}-\varpi_n'(t_k)
=&\,\int_{t_{k-1}}^{t_{k}}\frac{\varpi_n'(s)-\varpi_n'(t_k)}{\tau_{k}}\zd s
=-J_{n-k}^{(n)}
 \end{align*}
for $1\leq k\leq n-1$ $(2\leq n\leq N)$. Hence the desired first equality is obtained by a simple subtraction.
For the case of $k=n-1$ $(2\leq n\leq N)$, the above equality gives
\begin{align*}
a_{1}^{(n)}-\varpi_n'(t_{n-1})=-J_{1}^{(n)}.
 \end{align*}
We have $a^{(n)}_0=\frac{1-\theta}{1-\alpha}\varpi_n'(t_{n-1})$ such that
 $a^{(n)}_0-\varpi_n'(t_{n-1})=\frac{\theta}{1-2\theta}\varpi_n'(t_{n-1}).$
 Thus a simple subtraction yields the second equality and completes the proof.
 \end{proof}

\begin{lemma}\label{lem: ank diff bound}
If \Mss{1} holds, the positive coefficients $a^{(n)}_{n-k}$ in \eqref{eq: an} satisfy
 \begin{align*}
a_{n-k-1}^{(n)}-a_{n-k}^{(n)}
\geq\begin{cases}
	 b_{n-2}^{(n)}+\frac{6}{5}I_{n-1}^{(n)},
	&k=1,\vspace{0.1cm}\\
	 b_{n-k-1}^{(n)}+\rho_{k-1}b^{(n)}_{n-k+1}+\frac{1}{5}I_{n-k}^{(n)},
	&k>1,
\end{cases}
 \end{align*}
for $1\le k\le n-2$ $(3\leq n\leq N)$, and for $k=n-1$ $(2\leq n\leq N)$,
 \begin{align*}
a_{0}^{(n)}-a_{1}^{(n)}
\geq\begin{cases}
	 I_{1}^{(2)},
	&n=2,\vspace{0.1cm}\\
\rho_{n-2}b^{(n)}_{2}+I_{1}^{(n)},
	&n>2.
\end{cases}	
 \end{align*}
\end{lemma}
\begin{proof}
For fixed $n$, applying \cref{lem: bn qn2}~(i) and \cref{lem: qn qn1}~(i), we
obtain
\begin{align}\label{ank diff bound-1}
I_{n-k-1}^{(n)}
&=\frac{\rho_{k+1}I_{n-k-1}^{(n)}}{1+\rho_{k+1}}
	+\frac{I_{n-k-1}^{(n)}}{1+\rho_{k+1}}
\geq b_{n-k-1}^{(n)}+\frac{I_{n-k}^{(n)}}{\rho_{k}(1+\rho_{k+1})}\\
&\geq b_{n-k-1}^{(n)}+\frac{I_{n-k}^{(n)}}{\rho(1+\rho)}
=b_{n-k-1}^{(n)}+\frac{16}{77}\,I_{n-k}^{(n)},
\quad 1\leq k\leq n-2,\nonumber
 \end{align}
where the assumption \Mss{1} was used. By using \cref{lem: qn qn1}~(ii) and
\cref{lem: bn qn2}~(ii),
\begin{align*}
\frac{\rho_{k-1}^3}{2(1+\rho_{k-1})}J_{n-k}^{(n)}
\geq\frac{\rho_{k-1}^2}{2(1+\rho_{k-1})}J_{n-k+1}^{(n)}
\geq \rho_{k-1}b^{(n)}_{n-k+1},\quad 2\leq k\leq n-1.
 \end{align*}
Then, noting that $2+2x-x^3\geq 9/64$ for $x\in[0,7/4]$, we apply
\cref{lem: bn qn2}~(iii) and the assumption \Mss{1} to get
\begin{align}\label{ank diff bound-2}
J_{n-k}^{(n)}&=\frac{\rho_{k-1}^3}{2(1+\rho_{k-1})}J_{n-k}^{(n)}
	+\frac{2+2\rho_{k-1}-\rho_{k-1}^3}{2(1+\rho_{k-1})}J_{n-k}^{(n)}
	\geq\rho_{k-1}b^{(n)}_{n-k+1}+\frac{9}{352}I_{n-k}^{(n)},
 \end{align}
where $2\leq k\leq n-1.$ Hence, with help of
\eqref{ank diff bound-1}--\eqref{ank diff bound-2}, we apply
\cref{lem: ank diff} to find
\begin{align*}
a_{n-k-1}^{(n)}-a_{n-k}^{(n)}&=I_{n-k-1}^{(n)}+J_{n-k}^{(n)}
\geq b_{n-k-1}^{(n)}+\rho_{k-1}b_{n-k+1}
	+\braB{\frac{9}{352}+\frac{16}{77}}I_{n-k}^{(n)}\\
&>b_{n-k-1}^{(n)}+\rho_{k-1}b_{n-k+1}^{(n)}+\frac{1}{5}I_{n-k}^{(n)},
	\quad 2\leq
k\leq n-2.
 \end{align*}
If $k=1$, by applying \cref{lem: ank diff} with the bound \eqref{ank diff bound-1}
and \cref{lem: bn qn2} (iii), one has
\begin{align*}
a_{n-2}^{(n)}-a_{n-1}^{(n)}=&\,I_{n-2}^{(n)}+J_{n-1}^{(n)}
\geq b_{n-2}^{(n)}+\frac{16}{77}I_{n-1}^{(n)}+I_{n-1}^{(n)}\geq
b_{n-2}^{(n)}+\frac{6}{5}I_{n-1}^{(n)}.
 \end{align*}

To complete the proof, it remains to consider the case of $k=n-1$ $(2\leq n\leq N)$.
If $n=2$, \cref{lem: ank diff} and \cref{lem: bn qn2} (iii) yield
 \begin{align*}
a_{0}^{(2)}-a_{1}^{(2)}
=\frac{\theta}{1-2\theta}\varpi_n'(t_{1})+J_{1}^{(2)}>J_{1}^{(2)}>I_{1}^{(2)}.
 \end{align*}
Now treat the last case of $n\ge3$.
We apply \cref{lem: bn an} (by taking $k=n-2$), \cref{lem: an} (i) and
the given condition \Mss{1} to get
\begin{multline*}
\rho_{n-2}b_{2}^{(n)}
    \leq\frac{\theta\tau_{n-2}}{2(t_{n-\theta}-t_{n-2})}
        \frac{\rho_{n-2}^2}{1+\rho_{n-2}}a_{2}^{(n)}
    \leq\frac{\theta\rho_{n-2}}{2}
        \frac{\rho_{n-2}^2}{1+\rho_{n-2}}a_{2}^{(n)}
        \leq\frac{\theta\rho^3}{2(1+\rho)}a_{0}^{(n)}\\
	=\frac{343}{352}\,\theta a^{(n)}_0
    <\frac{\theta}{\tau_n}\omega_{2-\alpha}(t_{n-\theta}-t_{n-1})
	=\frac{\theta(1-\theta)}{1-2\theta}\varpi_n'(t_{n-1})
	\leq \frac{\theta}{1-2\theta}\varpi_n'(t_{n-1}).
\end{multline*}
Therefore  \cref{lem: ank diff} and \cref{lem: bn qn2} (iii) lead to
 \begin{align*}
a_{0}^{(n)}-a_{1}^{(n)}
=\frac{\theta}{1-2\theta}\varpi_n'(t_{n-1})+J_{1}^{(n)}>\rho_{n-2}b_{2}^{(n)}+I_{1}^{(n)}.
 \end{align*}
The proof is completed.
\end{proof}

Recalling the definition~\eqref{eq: weights}, we proceed to
apply \cref{lem: ank diff,lem: ank diff bound}.

\begin{proof}[Proof of \cref{thm: A1 A2} (II)]
With the notation $I_{n-k}^{(n)}$ defined in~\eqref{ank diff bound-notations},
we can write the desired inequality as
\[
A^{(n)}_{n-k-1}-A^{(n)}_{n-k}
    \ge(1+\rho_k)b^{(n)}_{n-k}+\frac1{5}I_{n-k}^{(n)},
    \quad1\le k\le n-1,
\]
and treat four separate cases covering all possibilities.
Indeed, from the definition~\eqref{eq: weights} of $A^{(n)}_{n-k}$, it is not
difficult to verify that
\begin{itemize}
  \item[(1)] If $k=1$ for $n=2$,
  $$A^{(2)}_0-A^{(2)}_1=(1+\rho_1)b^{(2)}_1+a^{(2)}_0-a^{(2)}_1;$$
  \item[(2)] If $k=n-1$ for $n\geq3$,
  $$A^{(n)}_0-A^{(n)}_1=(1+\rho_{n-1})b^{(n)}_1+a^{(2)}_0-a^{(2)}_1-\rho_{n-2}b^{(n)}_2;$$
  \item[(3)] If $k=1$ for $n\geq3$,
  $$A^{(n)}_{n-2}-A^{(n)}_{n-1}=(1+\rho_1)b^{(n)}_{n-1}+a^{(n)}_{n-2}-a^{(n)}_{n-1}-b^{(n)}_{n-2};$$
  \item[(4)] If $2\leq k\leq n-2$ for $n\geq4$,
  $$A^{(n)}_{n-k-1}-A^{(n)}_{n-k}
	=(1+\rho_k)b^{(n)}_{n-k}+a^{(n)}_{n-k-1}-a^{(n)}_{n-k}-b^{(n)}_{n-k-1}-\rho_{k-1}b^{(n)}_{n-k+1}\,.$$
\end{itemize}
The claimed inequality follows from \cref{lem: ank diff bound} directly
and completes the proof.
\end{proof}

\begin{proof}[Proof of \cref{thm: A1 A2} (III)]
The proof of \cref{lem: ank diff} shows that
\[
\frac{1-2\theta}{1-\theta}a^{(n)}_0-a^{(n)}_1=J^{(n)}_1>0
	\quad\text{for~$2\le n\le N$.}
\]
In the case~$n=2$, the definition~\eqref{eq: weights} gives
\begin{align*}
\frac{1-2\theta}{1-\theta}\,A^{(2)}_0-A^{(2)}_1
	&=\frac{1-2\theta}{1-\theta}\,\bigl(a^{(2)}_0+\rho_1b^{(2)}_1\bigr)
	-\bigl(a^{(2)}_1-b^{(2)}_1\bigr)\\
	&=J^{(2)}_1+\frac{1-2\theta}{1-\theta}\,\rho_1b^{(2)}_1+b^{(2)}_1>0.
\end{align*}
For~$n\ge3$, one has
\begin{align*}
\frac{1-2\theta}{1-\theta}\,A^{(n)}_0-A^{(n)}_1
	&=\frac{1-2\theta}{1-\theta}\,\bigl(a^{(n)}_0+\rho_{n-1}b^{(n)}_1\bigr)
	-\bigl(a^{(n)}_1+\rho_{n-2}b^{(n)}_2-b^{(n)}_1\bigr)\\
	&=J^{(n)}_1-\rho_{n-2}b^{(n)}_1
		+\frac{1-2\theta}{1-\theta}\,\rho_{n-1}b^{(n)}_1+b^{(n)}_1>0
\end{align*}
because $J^{(n)}_1\ge\rho_{n-1}b^{(n)}_2+\frac{9}{352}I^{(n)}_1$ from the
case~$k=n-1$ of~\eqref{ank diff bound-2}.
\end{proof}

\section{Numerical experiments}\label{sec: numerical}

A numerical example is reported here to support our theory.
The fully discrete scheme \eqref{eq: discrete IBVP} is used to
solve the subdiffusion problem \eqref{eq: IBVP}
in the domain $\Omega=(0,\pi)$ and $T=1$.
We take~$\kappa=2$ and set the exact solution
 $u(x,t)=\brab{1+\omega_{1+\sigma}(t)}\sin(x)$.
This solution satisfies a stronger estimate than~\eqref{eq: sigma}, namely,
$\mynorm{u^{(\nu)}(t)}_{H^2(\Omega)}\le C_ut^{\sigma-\nu}$ for $0<t\le T$ and $\nu\in\{1,2,3\}$.
As noted in \cite[Remark 7]{LiaoLiZhang:2018}, the  graded
mesh $t_n=T\bra{n/N}^{\gamma}$ satisfying \Mss{1}--\Mss{2},
is optimal in resolving the initial singularity.

In our computations, a linear finite element approximation is applied on a uniform mesh for $\Omega$ with $M$ nodes.
As done in an earlier paper~\cite{LiaoLiZhang:2018}, we split the interval $[0,T]$
into two parts~$[0,T_0]\cup[T_0,T]$. In first part~$[0,T_0]$ we used the smoothly graded
mesh~$t_n=(n/N_0)^\gamma T_0$ for~$0\le n\le N_0$, while a uniform mesh with step size~$\taumax$
is used in the second part~$[T_0,T]$.  For a given total
number~$N$ of time levels, we put
\[
T_0\defeq 2^{-\gamma}\quad\text{and}\quad
N_0\defeq\biggl\lceil\frac{\gamma N}{2^\gamma-1+\gamma}\biggr\rceil
\quad\text{so that}\quad
\tau\defeq\frac{T-T_0}{N-N_0}\ge\frac{\gamma T_0}{N_0}\ge\tau_{N_0}.
\]
To avoid problems with roundoff, the discrete coefficients $a^{(n)}_{n-k}$~and
$b^{(n)}_{n-k}$ from~\eqref{eq: an} and \eqref{eq: bn-anlternative}, respectively,
were computed using adaptive Gauss--Kronrod quadrature.

\begin{table}[!htp]
\caption{Numerical temporal accuracy for $\sigma=1+\alpha$ and $\gamma=1$.}
\label{tab: sigma=1+alpha}
\begin{center}
\renewcommand{\arraystretch}{1.1}
\begin{tabular}{c|cc|cc|cc}
&\multicolumn{2}{c|}{$\alpha=0.4$, $\sigma=1.4$}
&\multicolumn{2}{c|}{$\alpha=0.6$, $\sigma=1.6$}
&\multicolumn{2}{c}{$\alpha=0.8$, $\sigma=1.8$} \\
\cline{2-3}    \cline{4-5} \cline{6-7}
$N$&$e(N)$ & Order &$e(N)$ & Order &$e(N)$&Order\\
\hline
64    &2.78e-04 &--  &2.32e-04&--       &1.62e-04&-- \\
128  &7.24e-05 &1.94&5.97e-05&1.96&4.13e-05&1.97 \\
256  &1.87e-05 &1.95&1.52e-05&1.97&1.05e-05&1.97 \\
512  &4.74e-06 &1.97&3.85e-06&1.98&2.68e-06&1.98 \\
1024&1.59e-06 &1.58&9.72e-07&1.99&6.80e-07&1.98 \\
2048&5.61e-07 &1.50&2.45e-07&1.99&1.73e-07&1.97 \\
4096&2.01e-07 &1.48&6.06e-08&2.02&4.52e-08&1.94 \\
8192&5.83e-08 &1.46&1.23e-08&1.98&1.01e-08&1.83 \\
\hline
$\min\{\gamma\sigma,2\}$& &1.40&&1.60&&1.80
\end{tabular}
\end{center}
\end{table}

\begin{table}[!htp]
\caption{Numerical temporal accuracy for $\sigma=1.2$ and $\alpha=0.4$.}
\label{tab: sigma=1.2}
\begin{center}
\renewcommand{\arraystretch}{1.1}
\begin{tabular}{c|cc|cc|cc}
&\multicolumn{2}{c|}{$\gamma=1$}
&\multicolumn{2}{c|}{$\gamma=5/3=\gammaopt$}
&\multicolumn{2}{c}{$\gamma=2$} \\
\cline{2-3}    \cline{4-5} \cline{6-7}
$N$& $e(N)$ & Order &$e(N)$ & Order &$e(N)$&Order \\
\hline
64  &2.98e-04 &--  &1.29e-04&--        &2.12e-04&-- \\
128 &8.52e-05 &1.81&3.08e-05&2.07&5.07e-05&2.07 \\
256 &2.97e-05 &1.52&7.38e-06&2.06&1.24e-05&2.03\\
512 &1.18e-06 &1.33&1.77e-06&2.05&3.02e-06&2.04 \\
1024&4.81e-06 &1.30&4.21e-07&2.07&7.22e-07&2.06 \\
2048&1.98e-06 &1.27&9.25e-08&2.19&1.65e-07&2.12 \\
\hline
$\min\{\gamma\sigma,2\}$& &1.20&&2.00&&2.00
\end{tabular}
\end{center}
\end{table}

\begin{table}[!htp]
\caption{Numerical temporal accuracy for $\sigma=0.8$ and $\alpha=0.4$.}
\label{tab: sigma=0.8}
\begin{center}
\renewcommand{\arraystretch}{1.1}
\begin{tabular}{c|cc|cc|cc}
&\multicolumn{2}{c|}{$\gamma=2$}
&\multicolumn{2}{c|}{$\gamma=5/2=\gammaopt$}
&\multicolumn{2}{c}{$\gamma=3$} \\
\cline{2-3}    \cline{4-5} \cline{6-7}
$N$& $e(N)$ & Order &$e(N)$&Order &$e(N)$&Order\\
\hline
64  &3.52e-04&--  &5.28e-04&--  &5.04e-04&--\\
128 &8.17e-05&2.11&1.22e-04&2.11&1.17e-04&2.09 \\
256 &1.93e-05&2.08&2.93e-05&2.07&2.83e-05&2.06\\
512 &4.54e-06&2.08&7.02e-06&2.06&6.86e-06&2.05 \\
1024&1.08e-06&2.07&1.69e-06&2.05&1.68e-06&2.02 \\
2048&3.27e-07&1.73&4.29e-07&1.98&4.28e-07&1.97 \\
\hline
$\min\{\gamma\sigma,2\}$& &1.60&&2.00&&2.00 \\
\end{tabular}
\end{center}
\end{table}

\begin{table}[!htp]
\caption{Numerical temporal accuracy for $\sigma=0.4$ and $\alpha=0.4$.}
\label{tab: sigma=0.4}
\begin{center}
\renewcommand{\arraystretch}{1.1}
\begin{tabular}{c|cc|cc|cc}
&\multicolumn{2}{c|}{$\gamma=2$}
&\multicolumn{2}{c|}{$\gamma=5/2$}
&\multicolumn{2}{c}{$\gamma=5=\gammaopt$} \\
\cline{2-3}    \cline{4-5} \cline{6-7}
$N$& $e(N)$ & Order &$e(N)$&Order &$e(N)$&Order\\
\hline
64  &8.30e-03&--  &4.61e-03&--  &2.04e-03&--\\
128 &4.53e-03&0.87&2.23e-03&1.00&4.82e-04&2.08 \\
256 &2.56e-03&0.83&1.11e-03&1.01&1.22e-04&2.11\\
512 &1.45e-03&0.82&5.51e-04&1.00&2.66e-05&2.08\\
1024&8.25e-04&0.81&2.74e-04&1.01&6.40e-06&2.05 \\
2048&4.71e-04&0.81&1.37e-04&1.00&1.58e-06&2.02 \\
\hline
$\min\{\gamma\sigma,2\}$&&0.80&&1.00&&2.00 \\
\end{tabular}
\end{center}
\end{table}

Since the $O(h^2)$ behaviour of the spatial error is standard, we
fixed $M=10^4$ so that the temporal error dominates when~$N\le2,048$.  Thus,
from \cref{thm: convergence}, we expect the $L_\infty(L_2)$-error $e(N)\defeq\max_{1\le n\le N}\|u^n_h-u(t_n)\|$ to
behave like $O(\tau^{\min\{\gamma\sigma,2\}})$.  We tested the sharpness of this
prediction by four scenarios:
\begin{itemize}[itemindent=0.8cm]
\item[Table~\ref{tab: sigma=1+alpha}:] $\sigma=1+\alpha$ and $\gamma=1$ with
fractional orders $\alpha=0.4$, $0.6$~and $0.8$.
\item[Table~\ref{tab: sigma=1.2}:] $\sigma=1.2$~and $\alpha=0.4$ with mesh
parameters $\gamma=1$, $5/3$~and $2$.
\item[Table~\ref{tab: sigma=0.8}:] $\sigma=0.8$~and $\alpha=0.4$ with mesh
parameters $\gamma=2$, $5/2$~and $3$.
\item[Table~\ref{tab: sigma=0.4}:] $\sigma=0.4$~and $\alpha=0.4$ with mesh
parameters $\gamma=2$, $5/2$~and $5$.
\end{itemize}

The empirical order of convergence, listed as ``Order'' in the tables,
was computed in the usual way by supposing that $e(N)\approx C_u\tau^{q}$ and
evaluating the convergence rate $q\approx \log_2[e(N)/e(2N)]$.  The
optimal mesh parameter~$\gammaopt\defeq 2/\sigma$ is the smallest value
of~$\gamma$ for which we expect second-order convergence; for
$\gamma>\gammaopt$
we still expect second-order convergence but with a constant factor that grows
with~$\gamma$.  The convergence behaviour is always as expected, but it
is interesting to observe that, for larger values of~$\sigma$ (corresponding to
a less singular solution), the order can be close to~$2$ on the coarser grids.
In such cases, the predicted convergence order is not observed until the total number $N$ of time levels is
quite large.

\section{Concluding remarks}

An \textbf{ECS} analysis, including a reasonable \textbf{ECS} hypothesis
and a global consistency error, is proposed for investigating a class of
numerical approximations to the Caputo fractional derivative, employing
piecewise polynomial interpolation on general nonuniform time meshes. The
global consistency error bound \eqref{ieq: globalConsistencyBound}
reveals an interesting behavior: \emph{the global approximate error of numerical
Caputo formula would be ``local" despite being naturally
nonlocal}, by choosing the time mesh
according to the error equidistribution principle. The effectiveness of the \textbf{ECS} analysis is shown for the
familiar L1 formula in \cite{LiaoLiZhang:2018,LiaoYanZhang:2019}, and for the
higher-order Alikhanov formula in this paper. In the latter case, the
theoretical properties of the discrete convolution kernels (\cref{thm: A1 A2})
and a new interpolation error formula with integral remainder for quadratic
polynomials (\cref{lem:quadraticInterpolationError}) are crucial to obtaining a
useful \textbf{ECS} bound of the local truncation error.

As our answer to the problem formulated in \cref{sec: introduction}, the
fractional Crank--Nicolson time-stepping scheme \eqref{eq: discrete IBVP}
for the linear reaction-subdiffusion equation~\eqref{eq: IBVP} is stable
(\cref{thm: stability}) and convergent (\cref{thm: convergence}) on general
nonuniform grids satisfying a mild restriction ($\rho=7/4$) on the local time
step-size ratio. Consequently, some adaptive step-size control
with $\tau_{\mathrm{next}}/\tau_{\mathrm{current}}\ge4/7$
is permitted to resolve certain complex behaviors of the solution arising in
nonlinear time-fractional PDEs, but a sudden, drastic reduction of the step size
should be avoided to ensure stability.  For the linear case, a sharp
$L^2$-norm error estimate of the numerical solution is also presented
(\cref{thm: convergence}), demonstrating that a graded mesh can effectively
resolve the initial singularity.

The present results lead naturally to another question: can an
\textbf{ECS} analysis be applied to other high-order numerical Caputo formulas
such as the BDF2-like (L1-2) approximation
\cite{GaoSunZhang:2014,LiaoMcLeanZhang:2019,LiaoLyuVongZhao:2016,LvXu:2016}?
We plan to address this issue in further research.

\section*{Acknowledgements}
Hong-lin Liao thanks Prof.~Ying Zhao for her valuable discussions and fruitful
suggestions, and the hospitality of Beijing Computational Science
Research Center during the period of his visit.

\appendix
\section{Proof of \cref{lem:quadraticInterpolationError}}
\label{appendix:quadraticInterpolationError}

For fixed $n$ and $1\leq k\leq n-1$, let $\ell_{k,j}(t)$ ($j=k-1,k,k+1$) be the
basis functions of quadratic Lagrange interpolation $\Pi_{2,k}v$ at the
points $t_{k-1}$, $t_{k}$ and $t_{k+1}$.
Firstly, we will express the interpolation error
$\brab{\widetilde{\Pi_{2,k}}v}(t)=v(t)-\bra{\Pi_{2,k}v}(t)$ in an integral form.
To do so, recall two basic properties of basis functions, $\ell_{k,j}(t_l)=\delta_{jl}$ and
\begin{align}\label{BasisFunction-property1}
\sum_{j=k-1}^{k+1}\ell_{k,j}(t)(t_{j}-t)^{\nu}=\delta_{0\nu},\quad \nu\in\{0,1,2\}\,,
\end{align}
where $\delta_{jl}$ and $\delta_{0\nu}$ are Kronecker delta functions.
Now applying the Taylor's expansion with integral remainder, one has
\begin{align*}
v(t_{j})&=\sum_{m=0}^2\frac{v^{(m)}(t)}{m!}(t_{j}-t)^m
+\frac12\int_t^{t_{j}}(t_{j}-s)^2v'''(s)\zd{s}\,,\quad j\in\{k-1,k,k+1\}\,.
\end{align*}
Since $\bra{\Pi_{2,k}v}(t)=\sum_{j=k-1}^{k+1}\ell_{k,j}(t)v(t_{j})$,
a simple combination with the three weights (basis functions) $\ell_{k,j}(t)$ ($j=k-1,k,k+1$) gives
the interpolation error
\begin{align}\label{quadraticInterpolationError}
\brab{\widetilde{\Pi_{2,k}}v}(t)=\frac12\sum_{j=k-1}^{k+1}\int_{t_{j}}^t\ell_{k,j}(t)(t_{j}-s)^2v'''(s)\zd{s},\quad t_{k-1}\leq t\leq t_{k+1},
\end{align}
because the property \eqref{BasisFunction-property1} implies that
\begin{align*}
\sum_{j=k-1}^{k+1}\ell_{k,j}(t)\sum_{m=0}^2\frac{v^{(m)}(t)}{m!}(t_{j}-t)^m
=\sum_{m=0}^2\frac{v^{(m)}(t)}{m!}\sum_{j=k-1}^{k+1}\ell_{k,j}(t)(t_{j}-t)^m=v(t)\,.
\end{align*}
Furthermore, differentiating both sides of \eqref{quadraticInterpolationError}, one
applies \eqref{BasisFunction-property1} again to get
\begin{align}\label{CauptoInterpolationError-derivative}
\brab{\widetilde{\Pi_{2,k}}v}'(t)=&\,\frac12v'''(t)\sum_{j=k-1}^{k+1}\ell_{k,j}(t)(t_{j}-t)^2
+\frac12\sum_{j=k-1}^{k+1}\int_{t_{j}}^t\ell_{k,j}'(t)(t_{j}-s)^2v'''(s)\zd{s}\nonumber\\
=&\,\frac12\sum_{j=k-1}^{k+1}\int_{t_{j}}^t\ell_{k,j}'(t)(t_{j}-s)^2v'''(s)\zd{s}=\sum_{j=k-1}^{k+1}L_{j}(v)
\end{align}
for $t_{k-1}\leq t\leq t_{k+1}$, where
\[
L_{j}(v):=\frac12\int_{t_{j}}^t\ell_{k,j}'(t)(t_{j}-s)^2v'''(s)\zd{s}\,,\quad
j\in\{k-1,k,k+1\}\,.
\]
Secondly, we express the required integration error
$\int_{t_{k-1}}^{t_k}q'(t)\brab{\widetilde{\Pi_{2,k}}v}'(t)\zd{t}$
in terms of $\widetilde{\Pi_{1,k}}q$ by using
\eqref{CauptoInterpolationError-derivative}.
Since
$$\ell_{k,k+1}'(t)=\frac{2}{(\tau_{k+1}+\tau_{k})\tau_{k+1}}\bra{t-t_{k-1/2}},$$
\cref{lem: integralFormula} yields
\begin{align*}
\int_{t_{k-1}}^{t_k}\ell_{k,k+1}'(t)q'(t)\zd t=&\,
2\int_{t_{k-1}}^{t_k}\frac{\brat{t-t_{k-1/2}}q'(t)\zd t}{(\tau_{k+1}+\tau_{k})\tau_{k+1}}
=-2\int_{t_{k-1}}^{t_k}\frac{\brab{\widetilde{\Pi_{1,k}}q}(t)\zd t}{(\tau_{k+1}+\tau_{k})\tau_{k+1}}\,.
\end{align*}
Thus applying the formula for~$L_{k+1}(v)$, we
exchange the order of integration to find
\begin{align}\label{CauptoInterpolationError-object1}
\qquad&\,\int_{t_{k-1}}^{t_k}q'(t)L_{k+1}(v)\zd{t} =
\frac12\int_{t_{k-1}}^{t_k}\ell_{k,k+1}'(t)q'(t)\zd{t}\int_{t_{k+1}}^t(t_{k+1}-s)^2v'''(s)\zd{s}\\
&\,\hspace{2cm}=-\frac12\int_{t_{k-1}}^{t_k}\ell_{k,k+1}'(t)q'(t)\zd{t}\int_{t_{k}}^{t_{k+1}}(t_{k+1}-s)^2v'''(s)\zd{s}\nonumber\\
&\,\hspace{2.4cm}-\frac12\int_{t_{k-1}}^{t_k}\ell_{k,k+1}'(t)q'(t)\zd{t}\int_{t_{k}}^{t}(t_{k+1}-s)^2v'''(s)\zd{s}\nonumber\\
&\,\hspace{2cm}=\int_{t_{k}}^{t_{k+1}}(t_{k+1}-s)^2v'''(s)\zd{s}
\int_{t_{k-1}}^{t_k}\frac{\brab{\widetilde{\Pi_{1,k}}q}(t)\zd t}{(\tau_{k+1}+\tau_{k})\tau_{k+1}}\nonumber\\
&\,\hspace{2.4cm}-\frac12\int_{t_{k-1}}^{t_k}\kbra{(t_{k+1}-s)^2\int_{t_{k-1}}^{s}\ell_{k,k+1}'(t)q'(t)\zd{t}}v'''(s)\zd{s}\,.\nonumber
\end{align}
Similarly, it is easy to check the following equality
\begin{align*}
\ell_{k,k-1}'(t)=&\,\frac{2t-t_{k}-t_{k+1}}{(\tau_{k+1}+\tau_{k})\tau_{k}}
=\frac{2}{(\tau_{k+1}+\tau_{k})\tau_{k}}\brab{t-t_{k-1/2}}-\frac{1}{\tau_{k}}\,,
\end{align*}
and it follows that
\begin{align*}
\int_{t_{k-1}}^{t_k}\ell_{k,k-1}'(t)q'(t)\zd t=&\,
2\int_{t_{k-1}}^{t_k}\frac{\brat{t-t_{k-1/2}}q'(t)\zd t}{(\tau_{k+1}+\tau_{k})\tau_{k}}
-\frac1{\tau_k}\int_{t_{k-1}}^{t_k}q'(t)\zd t
\nonumber\\
=&\,-2\int_{t_{k-1}}^{t_k}\frac{\brab{\widetilde{\Pi_{1,k}}q}(t)\zd t}{(\tau_{k+1}+\tau_{k})\tau_{k}}
-\frac1{\tau_k}\int_{t_{k-1}}^{t_k}q'(t)\zd t\,.
\end{align*}
So applying the formula of $L_{k-1}(v)$, we
exchange the order of integration to find
\begin{align}\label{CauptoInterpolationError-object2}
&\int_{t_{k-1}}^{t_k}q'(t)L_{k-1}(v)\zd{t} =
\frac12\int_{t_{k-1}}^{t_k}q'(t)\int_{t_{k-1}}^t\ell_{k,k-1}'(t)(t_{k-1}-s)^2v'''(s)\zd{s}\zd{t}\\
&\hspace{1cm}=\frac12\int_{t_{k-1}}^{t_k}(t_{k-1}-s)^2v'''(s)\zd{s}\int_{s}^{t_{k}}\ell_{k,k-1}'(t)q'(t)\zd t\nonumber\\
&\hspace{1cm}=\frac12\int_{t_{k-1}}^{t_k}(t_{k-1}-s)^2v'''(s)\zd{s}\int_{t_{k-1}}^{t_{k}}\ell_{k,k-1}'(t)q'(t)\zd t\nonumber\\
&\hspace{1.3cm}-\frac12\int_{t_{k-1}}^{t_k}(t_{k-1}-s)^2v'''(s)\zd{s}\int_{t_{k-1}}^{s}\ell_{k,k-1}'(t)q'(t)\zd t\nonumber\\
&\hspace{1cm}=-\int_{t_{k-1}}^{t_k}(t_{k-1}-s)^2v'''(s)\zd{s}
\int_{t_{k-1}}^{t_k}\frac{\brab{\widetilde{\Pi_{1,k}}q}(t)\zd t}{(\tau_{k+1}+\tau_{k})\tau_{k}}\nonumber\\
&\hspace{1.3cm}-\frac12\!\int_{t_{k-1}}^{t_k}\!\kbra{\int_{t_{k-1}}^{t_k}\!\!\frac{q'(t)}{\tau_k}\zd t
+\int_{t_{k-1}}^{s}\!\!\!\ell_{k,k-1}'(t)q'(t)\zd t}(t_{k-1}-s)^2v'''(s)\zd{s}.\nonumber
\end{align}
For the remaining term involving $L_{k}(v)$, one has
\begin{align}\label{CauptoInterpolationError-object3}
\int_{t_{k-1}}^{t_k}q'(t)L_{k}(v)\zd{t}
=&\,\frac12\int_{t_{k-1}}^{t_k}q'(t)\int_{t_{k}}^t\ell_{k,k}'(t)(t_{k}-s)^2v'''(s)\zd{s}\zd{t}\nonumber\\
=&\,-\frac12\int_{t_{k-1}}^{t_k}\kbra{(t_{k}-s)^2\int_{t_{k-1}}^{s}\ell_{k,k}'(t)q'(t)\zd t}v'''(s)\zd{s}.
\end{align}
Then collecting the three equalities
\eqref{CauptoInterpolationError-object1}--%
\eqref{CauptoInterpolationError-object3},
one applies the formula \eqref{CauptoInterpolationError-derivative} to get
\begin{multline*}
\int_{t_{k-1}}^{t_k}q'(s)\brab{\widetilde{\Pi_{2,k}}v}'(s)\zd s
=\int_{t_{k}}^{t_{k+1}}(t_{k+1}-s)^2v'''(s)\zd{s}
\int_{t_{k-1}}^{t_k}\frac{\brab{\widetilde{\Pi_{1,k}}q}(t)\zd t}{(\tau_{k+1}
    +\tau_{k})\tau_{k+1}}\\
-\int_{t_{k-1}}^{t_k}(s-t_{k-1})^2v'''(s)\zd{s}
\int_{t_{k-1}}^{t_k}\frac{\brab{\widetilde{\Pi_{1,k}}q}(t)\zd t}{(\tau_{k+1}
    +\tau_{k})\tau_{k}}+\int_{t_{k-1}}^{t_{k}}\mathcal{K}_q(s)v'''(s)\zd{s},
\end{multline*}
where the integral kernel
\begin{align}\label{CauptoInterpolationError-kernelDef}
\mathcal{K}_q(s):=
-\frac{(t_{k-1}-s)^2}{2\tau_k}\int_{t_{k-1}}^{t_k}q'(t)\zd t
-\frac12\sum_{j=k-1}^{k+1}(t_{j}-s)^2\int_{t_{k-1}}^{s}\ell_{k,j}'(t)q'(t)\zd{t}\,.
\end{align}
Finally, to complete the proof, it remains to verify
\begin{align}\label{CauptoInterpolationError-kernel}
\mathcal{K}_q(s)=\int_{t_{k-1}}^{s}\brab{\widetilde{\Pi_{1,k}}q}(t)\zd{t}\,,\quad t_{k-1}\leq s\leq t_{k}.
\end{align}
Differentiating the identity
$(t-s)^2=\sum_{j=k-1}^{k+1}(t_j-s)^2\ell_{k,j}(t)$, we have
\begin{equation*}
\sum_{j=k-1}^{k+1}(t_j-s)^2\ell_{k,j}'(t)=2(t-s)\,,\quad t_{k-1}\leq s\leq t_{k+1}.
\end{equation*}
Thus it follows from \eqref{CauptoInterpolationError-kernelDef} that
\begin{align*}
\mathcal{K}_q(s)=
-\frac{(t_{k-1}-s)^2}{2\tau_k}\int_{t_{k-1}}^{t_k}q'(t)\zd t
-\int_{t_{k-1}}^{s}(t-s)q'(t)\zd{t}\,.
\end{align*}
We see that $\mathcal{K}_q(t_{k-1})=0$ and
\[
\mathcal{K}_q'(s)=
q(s)-q(t_{k-1})-\frac{q(t_{k})-q(t_{k-1})}{\tau_k}\bra{s-t_{k-1}}
=\brab{\widetilde{\Pi_{1,k}}q}(s),\quad t_{k-1}\leq s\leq t_{k},
\]
which leads to the desired result \eqref{CauptoInterpolationError-kernel} immediately since
$\mathcal{K}_q(s)=\int_{t_{k-1}}^{s}\mathcal{K}_q'(t)\zd{t}$.

\bibliographystyle{siamplain}

\end{document}